\documentclass{article}
\pdfoutput=1
\usepackage{latexsym, verbatim, amscd, amsmath, amssymb, amsthm}
\usepackage{pinlabel}
\usepackage{graphicx}
\usepackage{hyperref}

\hypersetup{
    pdftitle={Quasi-isometries Between Tubular Groups},    
    pdfauthor={Christopher H. Cashen},     
    pdfkeywords={Tubular Groups}, 
}

\theoremstyle{plain}
\newtheorem{theorem}{Theorem}[section]
\newtheorem{lemma}{Lemma}[section]
\newtheorem{proposition}{Proposition}[section]
\newtheorem{corollary}{Corollary}[section]
\newtheorem*{summarytheorem}{Theorem}
\newtheorem*{refthm}{Theorem}
\newtheorem*{maintheorem}{Main Theorem}

\theoremstyle{remark}
\newtheorem*{claim}{Claim}
\newtheorem*{note}{Note}
\newtheorem*{remark}{Remark}
\newtheorem*{remarks}{Remarks}
\newtheorem*{sketch}{Sketch}
\theoremstyle{definition}
\newtheorem{definition}{Definition}[section]
\newtheorem{example}{Example}[section]

\def\makeautorefname#1#2{\expandafter\def\csname#1autorefname\endcsname{#2}}
\let\fullref\autoref

\makeautorefname{theorem}{Theorem} 
\makeautorefname{lemma}{Lemma} 
\makeautorefname{proposition}{Proposition} 
\makeautorefname{corollary}{Corollary} 
\makeautorefname{definition}{Definition}
\makeautorefname{example}{Example}
\makeautorefname{subsection}{Section}
\makeautorefname{subsubsection}{Section}

\makeatletter \let\c@lemma=\c@theorem \makeatother
\makeatletter \let\c@proposition=\c@theorem \makeatother
\makeatletter \let\c@corollary=\c@theorem \makeatother
\makeatletter \let\c@definition=\c@theorem \makeatother
\makeatletter \let\c@example=\c@theorem \makeatother

\DeclareMathOperator{\GL}{GL}

\DeclareMathOperator{\PGL}{PGL}
\DeclareMathOperator{\Stab}{Stab}
\DeclareMathOperator{\Hom}{Hom}

\DeclareMathOperator{\slope}{slope}
\DeclareMathOperator{\gpcenter}{Z}
\DeclareMathOperator{\Err}{Err}
\DeclareMathOperator{\FIN}{FIN}

\def\from{\colon\thinspace}
\def\into{\hookrightarrow}
\newcommand{\case}[2]{\smallskip \noindent \emph{#1}:\newline #2 \medskip}
\newcounter{stepcount}
\newcommand{\step}{\item}
\newenvironment{steps}
{
\begin{list}{Step \arabic{stepcount}:}{\setlength{\itemindent}{12pt}\setlength{\leftmargin}{0pt}}\usecounter{stepcount}}
{\end{list}}
\def\qua{\;\;\,}

\begin{document}
\title{Quasi-isometries Between Tubular Groups}
%
\author{
Christopher H. Cashen\\
Department of Mathematics\\
University of Utah\\
Salt Lake City, UT 84112-0090 USA\\
\texttt{cashen@math.utah.edu}\\
\href{http://www.math.utah.edu/~cashen}{\url{http://www.math.utah.edu/~cashen}}
}
\maketitle
\begin{abstract}
We give a method of constructing maps between tubular groups inductively according to a finite set of strategies.
This map will be a quasi-isometry exactly when the set of strategies satisfies certain consistency criteria.
Conversely, if there exists a quasi-isometry between tubular groups, then there is a consistent set of strategies for building a quasi-isometry between them.

For two given tubular groups there are only finitely many candidate sets of strategies to consider, so it is possible in finite time to either produce a consistent set of strategies or decide that such a set does not exist.
Consequently, there is an algorithm that in finite time decides whether or not two tubular groups are quasi-isometric.
\end{abstract}


\section{Introduction}
One line of attack in Gromov's program to classify finitely generated groups up to quasi-isometry has been to study splittings of groups.
Given a splitting of a group $G$ into a graph of groups, one would like to know whether there is a similar splitting for a finitely generated group $H$ quasi-isometric to $G$, and what constraints such splittings impose upon quasi-isometries between the groups.

A result of this type applies when $G$ and $H$ are accessible groups.
Paposoglu and Whyte showed that $G$ and $H$ are quasi-isometric if and only if they have the same number of ends and if, in terminal splittings over finite subgroups, they have the same sets of quasi-isometry classes of one ended vertex groups \cite{PapWhy02}.
In terms of the Bass-Serre trees for such groups, this is a local restriction.

When edge groups are infinite there may be large scale restrictions as well.
Consider, for example, graphs of groups where all the local groups are infinite cyclic, such as the Baumslag-Solitar groups.
Here there are no obvious local restrictions; all the vertex spaces are quasi-isometric to each other.
However, crossing an edge space contributes a \emph{height change}.
Quasi-isometries between graphs of $\mathbb{Z}$'s must induce coarsely height preserving quasi-isometries of the corresponding Bass-Serre trees \cite{Why01}.

Perhaps the next most complicated situation to consider would be graphs of groups with infinite cyclic edge groups and vertex groups free abelian of rank two.
Such a group can be thought of as the fundamental group of a finite 2-complex consisting of a disjoint union of tori glued together by annuli.
Martin Bridson has described such groups as ``tubular''.

Mosher, Sageev, and Whyte prove splitting rigidity results for graphs of coarse Poincar\'{e} duality groups under some additional hypotheses \cite{MosSagWhy03}, \cite{MosSagWhy04}.
We define a class of tubular groups that satisfy these hypotheses.
Such splittings are quasi-isometrically rigid.
Furthermore, the patterns of attachment of edge groups to vertex groups are quasi-isometry invariants of the groups.
These edge patterns give local quasi-isometry restrictions, and height change gives a large scale restriction, similar to the graph of $\mathbb{Z}$'s example.

We will show that we can cover the Bass-Serre tree of a tubular group by a collection of infinite subtrees, called P--sets, that intersect pairwise in at most a single vertex.
We organize the P--sets into a tree of P--sets, which is essentially the nerve of the aforementioned cover.
The P--sets are similar to the Bass-Serre trees of Baumslag-Solitar groups, and we can use similar methods to construct quasi-isometries between P--set spaces.
To construct quasi-isometries of tubular groups we will need to coordinate the quasi-isometries on each P--set space so that they piece together in such a way as to satisfy the local edge-to-vertex pattern rigidity and global height change restrictions.

To accomplish this, we define \emph{sets of strategies} for inductively building quasi-isometries in a way that will satisfy the local restrictions given by edge pattern rigidity.
We give consistency criteria for such a set of strategies that determine whether the large scale restrictions on height change are satisfied.
\fullref{T:Algorithm} shows that there is an algorithm that in finite time will either produce a consistent set of strategies or decide that no such set exists.

We prove a number of results that can be summarized as:
\begin{summarytheorem}
  The following are equivalent:
  \begin{enumerate}
  \item The tubular groups $G_1$ and $G_2$ are quasi-isometric.
\item There exists an allowable isomorphism between the trees of P--sets of $G_1$ and $G_2$.
\item There exists a consistent set of strategies for $G_1$ and $G_2$.\label{Item:strategy}
  \end{enumerate}
\end{summarytheorem}

Combining these results with \fullref{T:Algorithm}, we get the main result of this paper:
\begin{maintheorem}
There is an algorithm that will take graph of groups decompositions of two tubular groups and in finite time decide whether or not the groups are quasi-isometric.  
\end{maintheorem}

\subsection{Previous Examples of Tubular Groups} 
Examples of tubular groups have appeared in the literature in various guises.

The fundamental group of the Torus Complexes of Croke and Kleiner \cite{CroKle00} is an example of a tubular group. 
These Torus Complexes provided examples of finite, non-positively curved, homeomorphic 2-complexes whose universal covers have non-homeomorphic ideal boundaries.

Right-Angled Artin groups whose defining graphs are trees of diameter at least three are tubular groups.
Such groups are also graph-manifold groups and have been shown to be quasi-isometric to each other by Behrstock and Neumann as a special case of their quasi-isometry classification of graph-manifold groups \cite{BehNeu06}.
This result is recovered as \fullref{Cor:RAAG}

Examples of tubular groups were used in work of Brady and Bridson \cite{BraBri00}, as well as Brady, Bridson, Forester, and Shankar \cite{BraBriFor09}, to show that the isoperimetric spectrum is dense in $[2,\infty)$ and that $\mathbb{Q}\cap[2,\infty)\subset IP$, respectively.
In the latter work these examples were termed ``snowflake groups''.

Wise's group \[W=\left<a,b,x,y\mid [a,b]=1, x^{-1}ax=(ab)^2, y^{-1}by=(ab)^2\right>\] is a tubular group that is non-Hopfian and CAT(0) \cite{Wis96}. 
It is not known if this group is automatic.
Resolving this question would either provide an example of a non-Hopfian automatic group or a non-automatic CAT(0) group.

\subsection{Acknowledgements}
This paper is based on a thesis \cite{Cas07} submitted in partial fulfillment of the requirements for the 
doctoral degree at the Graduate College of the University of Illinois at Chicago, under the supervision of Kevin Whyte. Additional work was conducted during MSRI's Geometric Group Theory Program.

\section{Preliminaries}
This section contains standard definitions and constructions.
Notation will most closely resemble that of Mosher, Sageev, and Whyte \cite{MosSagWhy04}.
\subsection{Coarse Geometry}
Let $(X,d_X)$ and $(Y,d_Y)$ be metric spaces, and let $A$, $B\subset X$.
The closed $r$-neighborhood of $A$ in $X$ is denoted $N_r(A)$.
The set $A$ is \emph{coarsely contained in} $B$, $A\stackrel{c}{\subset}B$, if $\exists r$ such that $A\subset N_r(B)$.
Subsets $A$ and $B$ are \emph{coarsely equivalent}, $A\stackrel{c}{=}B$, if $A\stackrel{c}{\subset}B$ and $B\stackrel{c}{\subset}A$.
A subset $A$ is \emph{coarsely dense} in $X$ if $A\stackrel{c}{=} X$.

A subspace $C$ of $X$ is a \emph{coarse intersection} of $A$ and $B$, $C=A\stackrel{c}{\cap}B$, if, for sufficiently large $r$, $C\stackrel{c}{=}N_r(A)\cap N_r(B)$.

If $A$ and $B$ are subspaces of $X$ then $A$ \emph{crosses} $B$ in $X$ if, for all sufficiently large $r$, there are at least two components $C_1$ and $C_2$ of $X\setminus N_r(B)$ such that for each $i$ and every $s>0$ there is a point $c\in A\cap C_i$ with $N_s(c)\subset C_i$.

A map $f\from X\to Y$ is \emph{$K$-bilipschitz}, for $K\geq 1$, if \[\frac{1}{K}d_X(x,y)\leq d_Y(f(x),f(y)) \leq Kd_X(x,y)\]  for all $x,y\in X$. 
The map $f$ is a \emph{$(K,C)$-quasi-isometric embedding}, for $K\geq 1,\, C\geq 0$, if \[\frac{1}{K}d_X(x,y)-C\leq d_Y(f(x),f(y)) \leq Kd_x(x,y)+C\]  for all $x,y\in X$. 
Furthermore, $f$ is a \emph{$(K,C)$-quasi-isometry} if it is a $(K,C)$-quasi-isometric embedding and the image is $C$-coarsely dense in $Y$.

Two maps $f,g\from X\to Y$ are \emph{bounded distance} from each other, $f\stackrel{c}{=}g$, if there is a $C\geq 0$ such that, $\forall x\in X$, $d_Y(f(x),g(x))\leq C$. 
Two maps, $f\from X\to Y$ and $\bar{f}\from Y\to X$, are \emph{coarse inverses} if $f\circ \bar{f} \stackrel{c}{=}Id_Y$ and $\bar{f}\circ f \stackrel{c}{=}Id_X$.
If $f$ is a quasi-isometry, there is a coarse inverse $\bar{f}\from Y\to X$ of $f$ that is also a quasi-isometry, with constants depending on those of $f$.

\subsection{Bass-Serre Theory}\label{S:BSTheory}
If $\Gamma$ is a graph, let $\mathcal{V}\varGamma$ denote the vertex set, and $\mathcal{E}\varGamma$ the edge set. Let $\mathcal{VE}\varGamma = \mathcal{V}\varGamma \cup \mathcal{E}\varGamma$.
The set of endpoints of $\Gamma$ is $\mathcal{E}\varGamma \times \{0,1\}$.

Each edge has two endpoints, of the form $\eta=(e,i)$, where the $i$ should be taken mod 2.
Each endpoint is identified with some vertex $v(\eta)$ such that $e=e(\eta)$ is incident to $v(\eta)$.

A \emph{graph of groups}, $(\Gamma, \{G_{\gamma}\}, \{\phi_{\eta}\})$, is a graph, $\Gamma$, equipped with a \emph{local group} $G_{\gamma}$ for each $\gamma \in \mathcal{VE}\varGamma$, and \emph{edge injections} $\phi_{\eta}\in \Hom(G_{e(\eta)}, G_{v(\eta)})$ for each endpoint $\eta$.
We will generally use $\Gamma$ to denote the graph of groups, and refer to the \emph{underlying graph of $\Gamma$} if we wish to consider only the graph itself.

\begin{note}
A graph of groups is \emph{of finite type} if the underlying graph is finite, the vertex groups are finitely presented, and the edge groups are finitely generated.
All the graphs of groups of interest in this paper are of finite type, so, from this point forward, finite type can be taken as an implicit hypothesis for any statement about graphs of groups.
\end{note}

Associated to a graph of groups there is a finitely presented group, 
\[G=\pi_{1}(\Gamma, \{G_{\gamma}\},\{\phi_{\eta}\}),\] 
the \emph{fundamental group of the graph of groups} \cite{Ser03}.

A graph of groups is \emph{reducible} if there is an edge $e$ such that the vertices $v(e,0)$ and $v(e,1)$ are distinct, and such that one of the edge homomorphisms $\phi_{(e,i)}$ is surjective.
In this case it is possible to simplify the graph of groups without changing the fundamental group.
Remove $e$ and $v(e,i)$ from $\Gamma$, and for any other endpoint with $v(e',j)=v(e,i)$, replace $\phi_{(e',j)}$ with:
 \[\phi_{(e,i+1)}\circ \phi^{-1}_{(e,i)}\circ \phi_{(e',j)}\]

Scott and Wall gave a topological realization of $G$ \cite{ScoWal79}.
Build a \emph{graph of spaces}, $\mathcal{K}$, for $(\Gamma, \{G_{\gamma}\}, \{\phi_{\eta}\})$ by choosing \emph{local spaces}, $\mathcal{K}_{\gamma}$, for each $\gamma \in \mathcal{VE}\varGamma$.
For each $\gamma$, choose $\mathcal{K}_{\gamma}$ to be a pointed, connected, compact CW-complex, with a map $\pi_1(\mathcal{K}_{\gamma})\to G_{\gamma}$. 
This map should be an isomorphism if $\gamma \in \mathcal{V}\varGamma$, and an epimorphism if $\gamma \in \mathcal{E}\varGamma$.
For each endpoint $\eta$ of $\Gamma$, choose an \emph{edge map}, a pointed CW-map $f_{\eta}\from\mathcal{K}_{e(\eta)} \to \mathcal{K}_{v(\eta)}$, such that the induced map on fundamental groups agrees with the edge injection.

Now, let $\mathcal{K}$ be the finite CW-complex obtained from the disjoint union
\[\coprod_{v\in \mathcal{V}\varGamma} \mathcal{K}_v \amalg \coprod_{e\in \mathcal{E}\varGamma} \mathcal{K}_e\times [0,1]\]
by using $f_{(e,i)}$ to glue $\mathcal{K}_{e\times \{i\}}$ to $\mathcal{K}_{v(e,i)}$ for each endpoint $(e,i)$.
The fundamental group $\pi_1(\mathcal{K})$ is well defined, up to isomorphism, and, by van Kampen's Theorem, is isomorphic to $G$. 

Consider the universal cover $X=\tilde{\mathcal{K}}$, with covering map $p$ and metric lifted from $\mathcal{K}$.
The group $G$ acts properly discontinuously, cocompactly and isometrically by deck transformations on $X$, so $G$ and $X$ are quasi-isometric by the \u{S}varc-Milnor Lemma \cite{BriHae99}.
Thus, $X$ serves as a \emph{geometric model} for $G$.
For questions of the coarse geometry of $G$, it is sufficient to study $X$.

The space $X$ can be decomposed into path connected lifts of the local spaces $\mathcal{K}_\gamma$, and the action of $G$ on $X$ respects this decomposition.
Let $D\varGamma=q(X)$ be the quotient space of the decomposition.
The quotient is a tree on which $G$ acts without edge inversion, and $D\varGamma$ is $G$-equivariantly isomorphic to the Bass-Serre Tree of $\Gamma$.
We use the notation $D\varGamma$ because the tree is the ``development'' of $\Gamma$.
Call $X$ the \emph{Bass-Serre Complex}, and $X\stackrel{q}{\to} D\varGamma$ the \emph{Bass-Serre tree of spaces} for $\Gamma$.

For $v\in \mathcal{V}D\varGamma$, $X_v=q^{-1}(v)$ is called a \emph{vertex space}, and $\mathcal{V}X=\bigcup_{v\in \mathcal{V}D\varGamma} X_v$ is the set of vertex spaces.
The set of vertex spaces is $\frac{1}{2}$-coarsely dense in $X$.
For $e\in \mathcal{E}D\varGamma$, $X_e=q^{-1}(\text{midpoint of } e)$ is called an \emph{edge space}.

For $t\in D\varGamma$, $\Stab_G(t)=\Stab_G(X_t)$.
The group $\Stab_G(t)$ is conjugate in $G$ to $G_{p(t)}$, and $\Stab_G(t)$ acts on $X_t$ as the deck transformation group of the covering map $X_t\to \mathcal{K}_{p(t)}$.

\section{Tubular Groups}

\subsection{Definition and Rigidity Results}
\subsubsection{Definition of Tubular Group}
The motivating examples of tubular groups are graphs of groups with edge groups $\mathbb{Z}$ and vertex groups $\mathbb{Z}^2$.
However, this description is not sufficient to give a quasi-isometrically closed class of groups.

\begin{definition}\label{Definition:tubular}
A \emph{tubular group} is the fundamental group of a finite, connected graph of groups satisfying the following conditions:
\begin{enumerate}
\item Every edge group is finitely generated and quasi-isometric to $\mathbb{Z}$.
\item Every vertex group is finitely generated and quasi-isometric to either $\mathbb{Z}$ or $\mathbb{Z}^2$.
\item There is at least one vertex group quasi-isometric to $\mathbb{Z}^2$ and at least one edge.\label{Li:reallytubular}
\item Every vertex whose local group is quasi-isometric to $\mathbb{Z}^2$ satisfies the crossing graph condition. (see below)\label{Li:incidentedgescross}
\end{enumerate}
\end{definition}
\begin{remarks}\hfill
\begin{itemize}
\item Graphs of groups with all local groups quasi-isometric to $\mathbb{Z}$ have been classified by Whyte \cite{Why01}. Conditions (\ref{Li:reallytubular}) and (\ref{Li:incidentedgescross}) exclude such groups from consideration.

\item Suppose $v\in \mathcal{V}\varGamma$ with $G_v$ quasi-isometric to $\mathbb{Z}$. Suppose $e\in \mathcal{E}\varGamma$ such that $v(e,i)=v$ and $v(e,i+1)\neq v$.
We can assume that $\phi_{(e,i)}(G_e)$ is a subgroup of $G_v$ of index at least two.
Otherwise, the graph of groups is reducible and could be simplified without changing the fundamental group.

\end{itemize}
\end{remarks}
\subsubsection{The Crossing Graph Condition}
In a graph of groups, a \emph{depth zero} vertex group is one that is not strictly coarsely contained in any other vertex group.

In a tubular group, those are precisely the vertex groups quasi-isometric to $\mathbb{Z}^2$.

Let $\mathcal{V}_0D\varGamma$ be the set of depth zero vertices of $D\varGamma$.
Let $\mathcal{V}_0X=\bigcup_{v\in \mathcal{V}_0D\varGamma}X_v$.

We define a \emph{crossing graph condition} for tubular groups.
Mosher, Sageev, and Whyte define a crossing graph condition in \cite{MosSagWhy04} that applies to more general graphs of groups. 
For tubular groups the two definitions are equivalent.

\begin{definition}
For each $v\in\mathcal{V}_0D\varGamma$ the \emph{crossing graph} for $v$ is a graph with one vertex for each edge incident to $v$ in $D\varGamma$.
Vertices of the crossing graph corresponding to edges $e$ and $e'$ are joined by an edge in the crossing graph if either $X_e$ and $X_{e'}$ cross in $X_v$ or if there is a third edge $e''$ incident to $v$ such that $X_{e''}$ crosses both $X_e$ and $X_{e'}$ in $X_v$. 
\end{definition}

\begin{definition}
  A depth zero vertex satisfies the \emph{crossing graph condition} if its crossing graph is connected.
\end{definition}

If the vertex group is isomorphic to $\mathbb{Z}^2$ the crossing graph condition is even simpler: the vertex satisfies the crossing graph condition if and only if the incident edge groups rationally span the vertex group.

The crossing graph condition fails in a tubular group only if for some depth zero vertex $v$, all of the incident edges have edge spaces coarsely equivalent to each other in $X_v$.

For example, the group $F_2\times \mathbb{Z}=\left<a,b\right>\times \left<c\right>$ can be regarded as a graph of groups with one $\mathbb{Z}^2$ vertex and one $\mathbb{Z}$ edge in two different ways. 
The $\mathbb{Z}^2$ vertex group could be either $\left<a,c\right>$ or $\left<b,c\right>$. 
Neither of these descriptions satisfy the crossing graph condition, because in each case the images of the two edge injections are the same cyclic subgroup of the vertex group.
The quasi-isometry coming from the obvious isomorphism does not respect vertex spaces, it takes the vertex space corresponding to $\left<a,c\right>$ to a union of edge spaces.

The crossing graph condition will prevent this sort of problem; it will ensure that a quasi-isometry of tubular groups takes depth zero vertex spaces to within bounded distance of depth zero vertex spaces. We will not give a direct proof of this statement for tubular groups, it will be a consequence of more general theorems of Mosher, Sageev, and Whyte. However, it is easy to get an idea of why this works. 
If a vertex satisfies the crossing graph condition then no line, $l$, in the vertex space can coarsely separate the vertex space in the whole tree of spaces. There will be some transverse line with an edge space attached, so we could avoid any neighborhood of $l$ by going far out in either direction along the transverse edge line, then crossing to the other side of the edge strip. On the other side of the edge strip is a vertex space quasi-isometric to a plane, and the intersection of any neighborhood of $l$ with this plane is bounded, so we can avoid it in the plane. If the image of a depth zero vertex space is not contained in a bounded neighborhood of some depth zero vertex space then it must cross deeply into components on opposite sides of some edge space. That edge space coarsely separates the entire tree of spaces, which means its preimage would be a line coarsely separating the vertex space in the tree of spaces, and this would be a contradiction.
\subsubsection{Rigidity Results of Mosher-Sageev-Whyte}\label{S:rigidity}
In this subsection we recall results of Mosher, Sageev, and Whyte \cite{MosSagWhy04} and apply them to tubular groups.

Given a graph of groups, $\Gamma$, with Bass-Serre tree of spaces $X \to T$, the following hypotheses will be needed: 
\begin{enumerate}
\item $\Gamma$ is finite type, irreducible, and finite depth.
\item No depth zero raft of the Bass-Serre tree $T$ is a line.
\item Each depth zero vertex group is coarse PD.
\item The crossing graph condition holds for each depth zero vertex of $T$ that is a raft. \label{Li:crossinggraph}
\item Each vertex and edge group of $\Gamma$ is coarse finite type.
\end{enumerate}

\begin{proposition}
Tubular groups satisfy these hypotheses.
\end{proposition}
\begin{proof}
The depth zero vertices are all rafts, and these are the only depth zero rafts.
We have included the crossing graph condition in the definition of tubular groups.
We can assume irreducibility, as discussed in the remarks following the definition of tubular groups.
Virtually abelian groups are coarse PD and coarse finite type.
Finite depth means there is a bound to the length of a chain of proper coarse inclusions of vertex and edge spaces. 
For tubular groups there are only chains of length two.
\end{proof}

\begin{theorem}[Quasi-isometric Rigidity Theorem]{\rm \cite[Theorem 1.5]{MosSagWhy04}}\qua
Let $\Gamma$ be a graph of groups satisfying (1)-(5) above. 
If $H$ is a finitely generated group quasi-isometric to $\pi_1 \Gamma$ then $H$ is the fundamental group of a graph of groups satisfying (1)-(5).
\end{theorem}

\begin{theorem}[Quasi-isometric Classification Theorem]{\rm \cite[Theorem 1.6]{MosSagWhy04}}\label{T:MSWrigidity}\qua
Let $\Gamma$, $\Gamma'$ be graphs of groups satisfying (1)-(5) above. 
Let $X\to T$, $X'\to T'$ be Bass-Serre trees of spaces for $\Gamma$, $\Gamma'$, respectively. 
If $f\from X\to X'$ is a quasi-isometry then $f$ coarsely respects vertex and edge spaces. 
To be precise, for any $K\geq 1$, $C\geq 0$ there exists a $K'$, $C'$ quasi-isometry $f_{\#}\from\mathcal{VE}(T)\to \mathcal{VE}(T')$ such that the following hold:
\begin{itemize}
\item If $a\in \mathcal{VE}(T)$ then $d_{\mathcal{H}}(f(X_a),X'_{f_{\#}(a)})\leq C'$
\item If $a'\in \mathcal{VE}(T')$ then there exists $a\in \mathcal{VE}(T)$ such that $d_{\mathcal{H}}(f(X_a),X'_{a'})\leq C'$
\end{itemize}
\end{theorem}

\begin{corollary}\label{corollary:classoftubulargroupsisclosed}
The class of tubular groups is closed under quasi-isometry. That is, any finitely generated group quasi-isometric to a tubular group is itself a tubular group.
Furthermore, any quasi-isometry between tubular groups coarsely respects vertex and edge spaces.
\end{corollary}

Note that in a tubular group, a vertex space quasi-isometric to $\mathbb{Z}^2$ is not bounded Hausdorff distance from any other vertex space.
Thus, we can change a quasi-isometry by a bounded amount so that it actually respects such vertex spaces. 

If $V\subset \mathbb{R}^n$ is a linear subspace, let $P_V$ be the set of affine subspaces of $\mathbb{R}^n$ parallel to $V$.
For a finite collection, $F$, of linear subspaces, the \emph{affine pattern} induced by $F$, $P_F$, is the union of the $P_V$ for $V\in F$.
An affine pattern is \emph{rigid} if for every $K,C,R$ there is an $R'$ such that if $f\from\mathbb{R}^n\to \mathbb{R}^n$ is a $K,C$ quasi-isometry that $R$-coarsely respects each $P_V$, then $f$ is within $R'$ of an affine homothety.

\begin{lemma}{\rm \cite[Lemma~7.2]{MosSagWhy04}}\qua
If $F$ is a finite collection of linear subspaces of $\mathbb{R}^n$ that contains $n+1$ hyperplanes in general position, then $F$ is rigid.
\end{lemma}

In particular, a collection of at least three distinct lines in the plane is rigid.

\begin{corollary}{\rm \cite[Corollary 7.11]{MosSagWhy04}}\label{Corollary:patterninvariance}\qua
Let $\Gamma$ be a graph of groups with all vertex and edge groups finitely generated abelian groups. 
Assume that for each depth zero, one vertex raft $v$ in the Bass-Serre tree, the collection of edge spaces at the vertex space of $v$ is a rigid affine pattern. 
Assume also that there are no line-like rafts of depth zero. If $H$ is any finitely generated group quasi-isometric to $G=\pi_1(\Gamma)$, then $H$ splits as a graph of virtually abelian groups and the quasi-isometry $G\to H$ is affine along each depth zero, one vertex raft.
Moreover, the set of affine equivalence classes of edge patterns is the same for $H$ as it is for $G$.
\end{corollary}

\subsection{Geometric Models for Tubular Groups}\label{S:geometricmodel}
Recall that an edge pattern is rigid if any quasi-isometry that respects each of the families of parallel lines is bounded distance from an affine homothety.
We will choose the metrics on the vertex spaces so that rigid patterns are ``symmetric''. 
The payoff for choosing the metrics in this way will be that any quasi-isometry which respects the edge pattern (possibly permuting families of parallel lines) is bounded distance from the composition of an isometry and a homothety. As a consequence, a quasi-isometry of tubular groups restricted to a depth zero vertex space with an edge pattern consisting of at least three families of parallel lines must stretch distances by the same amount in every direction.

\subsubsection{Affine Patterns of Lines in the Plane}
Let $F=\{l_1,\,l_2,..,\,l_n\}$, $n\geq 3$, be a finite collection of distinct lines through the origin in $\mathbb{R}^2$, with the usual Euclidean metric.
Two such collections, $F$ and $F'$, are \emph{linearly equivalent} if there exists $A\in \GL_2\mathbb{R}$ such that $AF=\{Al_i\}=F'$.
Scalar matrices in $\GL_2\mathbb{R}$ will preserve any such $F$, so we projectivise.

Let $P F=\{m_1,\,m_2,..,\, m_n\}$ be a collection of slopes in $\mathbb{RP}$.
There is a finite subgroup, $L_{P F}\subset \PGL_2\mathbb{R}$, that fixes $P F$ set-wise.

Lift $L_{PF}$ to $\GL_2(\mathbb{R})$ by taking $L_F=L_{P F}\gpcenter (\GL_2\mathbb{R})\cap \{A\in \GL_2\mathbb{R}\mid \det(A)=\pm 1\}$. 
We call $L_{P F}$ the \emph{group of symmetries} of $F$, and call $F$ \emph{symmetric} if $L_F$ acts by isometries.

\begin{proposition}\label{P:choosesymmetric}
For any collection $F$ of $n\geq 3$ distinct lines through the origin, there is a symmetric representative of the linear equivalence class of $F$. 
Choosing an isometry class of symmetric representative is equivalent to choosing a Euclidean metric on $\mathbb{R}^2$ for which $F$ is symmetric.
\end{proposition}
\begin{proof}
$L_F$ is finite.
Define a new metric on the plane by \[\left<x,y\right>_{L_F}=\frac{1}{|L_F|}\sum_{A\in L_F} \left<Ax, Ay\right>.\]
The group $L_F$ acts isometrically on the plane with this metric.
There is an $A_F\in \GL_2\mathbb{R}$ such that $\left<x,y\right>_{L_F}=\left<A_Fx,A_Fy\right>$, so $A_FF$ is a symmetric representative of the linear equivalence class of $F$.

Conversely, suppose $F'$ is a symmetric representative for the linear equivalence class of $F$.
Suppose $A_1$ and $A_2$ are matrices such that $|\det(A_i)|= 1$ and $A_iF=F'$.
Then $A_1A_2^{-1}\in L_{F'}$.
Since $F'$ is symmetric, this means the $A_i$ differ by an isometry, so 
\[\left< A_1x,A_1y\right> =\left< A_2x,A_2y\right>.\qedhere\] 
\end{proof}
For $n=2$ we will consider two orthogonal lines to be a symmetric pattern, but the group of symmetries in this case is not a finite group.

For $n=3$ there is a single linear equivalence class. 
There is, up to isometry, a unique symmetric representative, consisting of three lines meeting at angles $\frac{\pi}{3}$.
The group of symmetries is isomorphic to $S_3$.

For $n=4$ there are infinitely many linear equivalence classes, indexed by the cross ratios $[-1,0)$.
Each has, up to isometry, a unique symmetric representative, and a transitive group of symmetries.
When the cross ratio is $-1$ the group of symmetries is isomorphic to the dihedral group of order 8.
Otherwise, the group of symmetries is isomorphic to $\mathbb{Z}/2\mathbb{Z}\times \mathbb{Z}/2\mathbb{Z}$.

For $n=5$ there is no longer a unique symmetric representative for every linear equivalence class.
Indeed, there are five-line patterns with trivial group of symmetries.
In such a class, every member is symmetric, so the pattern does not determine a canonical choice of metric on the plane.

\subsubsection{Coarse Bass-Serre Complex}
It is possible to relax some of the group theoretic restrictions from \fullref{S:BSTheory} and still get a geometric model quasi-isometric to $G=\pi_1\Gamma$.
Following the proof of Lemma~2.9 of \cite{MosSagWhy04}, we construct a \emph{coarse Bass-Serre complex} $X'\to D\varGamma$ to serve as a geometric model for $G$.

Let $\Gamma$ be a graph of groups for a tubular group $G=\pi_1\Gamma$.

Let $v$ be a depth zero vertex of $\Gamma$, and $e$ an edge of $\Gamma$ incident to $v$ at endpoint $\eta$.

Suppose $A$ is a group quasi-isometric to $\mathbb{Z}^n$.
The group $A$ has a finite index normal subgroup $B\cong \mathbb{Z}^n$ \cite{BriGer96}.
There is a quasi-isometry, $f_B\from A\to B$, at bounded distance from $Id_A$, with $f_B|_B=Id_B$.
Furthermore, $f_B$ takes cyclic subgroups of $A$ to within bounded distance of cyclic subgroups of $B$.

Applying this reasoning to $G_v$ and $G_e$, 
\[\begin{CD}
G_e @>f_{\mathbb{Z}}>qi> \mathbb{Z}\\
@V\phi_\eta VV @VVV\\
G_v @>f_{\mathbb{Z}^2}>qi> \mathbb{Z}^2
\end{CD}\]

the image of the map $f_{\mathbb{Z}^2}\circ \phi_\eta \circ \bar{f}_{\mathbb{Z}}$ is bounded distance from a cyclic subgroup 
\[\left<x^ay^b\right><\left<x,y\right> \cong \mathbb{Z}^2.\]

The usual inclusion of $\mathbb{Z}^2$ into $\mathbb{R}^2$ is a quasi-isometry, and the subgroup $\left<x^ay^b\right>$ includes into the line through the origin with rational slope $\frac{b}{a}$.
In this way, each edge incident to $v$ is associated to a line in $\mathbb{R}^2$.
Let $F_v$ be the set of distinct lines.
The affine pattern induced by this set is called the \emph{edge pattern}, and is well defined up to linear equivalence.
If $F_v$ contains $n$ distinct lines, we will say that $v$ \emph{has n lines} or \emph{is an n-line vertex}.
Furthermore, we can choose a new Euclidean metric on $\mathbb{R}^2$ to make this pattern symmetric.

Let $X\stackrel{q}{\to}D\varGamma$ be a Bass-Serre tree of spaces for $\Gamma$.
For each $\gamma \in \mathcal{VE}D\varGamma$, let $h_\gamma \from X_{\gamma} \to X'_{{\gamma}}$, where $X'_{\gamma}$ is either $\mathbb{R}$ or $\mathbb{R}^2$ and $h_{\gamma}$ is the quasi-isometry given by restricting to a finite index normal abelian subgroup and then including into $X'_\gamma$.

For each vertex $v\in \mathcal{V}D\varGamma$, and each incident edge $e\in \mathcal{E}D\varGamma$, let $F_{ev} \from X_e\to X_v$ be the attaching map.
Define a new attaching map by $F'_{ev}=h_v\circ F_{ev}\circ \bar{h}_e \from X'_e \to X'_v$.

Define $X'$ by taking the disjoint union
\[\coprod_{v\in \mathcal{V}D\varGamma} X'_v \amalg \coprod_{e\in \mathcal{E}D\varGamma} X'_e\times [0,1]\]
and gluing according to the new attaching maps.

The $h_\gamma$ have uniform quasi-isometry constants, since the underlying graph was finite, so they piece together to give a quasi-isometry $h\from X\to X'$.
Since $G$ was quasi-isometric to $X$, we now have that $G$ is quasi-isometric to $X'$.
The action of $G$ on $X$ is quasi-conjugated by $h$ to give a proper, cobounded quasi-action of $G$ on $X'$.
The space $X'$ is still a tree of spaces over $D\varGamma$, and $X'\to D\varGamma$ is called a coarse Bass-Serre complex.

When $v$ is a depth zero vertex, the edge spaces of the incident edges attach to $X'_v$ along lines of the edge pattern in $X'_v$.

\subsubsection{Contraction Factors and Height Change}
Let $e$ be an edge of $D\varGamma$, and let $v_i$ be the vertex at endpoint $(e,i)$ of $e$, for $i=0,1$.
For each $i$, $F'_{ev_i}$ maps $X'_e=\mathbb{R}$ within bounded distance of a line in $X'_{v_i}$, and there is a factor $l_i$ such that $d_{X'_{v_i}}(F'_{ev_i}(x),F'_{ev_i}(y))=l_id_{X'_e}(x,y)$ to within bounded additive error.
Define the \emph{contraction factor across $e$} to be $\frac{l_1}{l_0}$ and the \emph{height change across $e$} to be: 
\[h(e)=-\log_2(\frac{l_1}{l_0})\]

Metrize the strip $X'_e\times [0,1]$ as the strip $0\leq y\leq 1$ in the plane with metric $(\frac{l_1}{l_0})^{2y}dx^2+dy^2$, so the edge strips are horostrips in a plane of constant curvature $-(\ln \frac{1_1}{l_0})^2$. 
Without changing the quasi-isometry type, we can change $X'$ by a bounded amount so that $F'_{ev_i}$ glues $X'_e\times \{i\}$ isometrically along a line in $X'_{v_i}$.

Let $p$ be vertical projection from the bottom ($i=0$) of the strip to the top ($i=1$) of the strip.
The map $p$ is closest point projection.
If $a$ and $b$ are points in the bottom edge of the strip, and $d_i$ is the distance in $X'_e\times \{i\}$, then:
 \[\frac{d_1(p(a),p(b))}{d_0(a,b)}=\frac{l_2}{l_1}=2^{-h(e)}\]

For vertices $v,w\in D\varGamma$, the \emph{height change from $v$ to $w$}, $h(v,w)$, is the sum of the height changes across the edges of the geodesic between $v$ and $w$.
This quantity will sometimes be called the \emph{height of $w$ relative to $v$}.
\begin{remark}
  If an edge has height change $h$ that means that closest point projection across the edge scales distances by the same amount as closest point projection between two horizontal horocycles whose $y$-coordinates (heights) differ by $h$ in the plane with metric $(\frac{1}{2})^{2y}dx^2+dy^2$. 
The obvious choice would have been to define height change using the natural logarithm, in which case the analogy would be to the hyperbolic plane. 
The choice of $\log_2$ turns out to be the convenient choice when relating height change to isoperimetric function, but the theory is the same regardless of which base is chosen. 
\end{remark}
\begin{example}
Suppose we have an edge $e$ going from vertex $v$ to vertex $w$.
Suppose we have chosen metrics on $X_v$ and $X_w$ so that the stabilizer subgroups are the usual integer lattice in the plane, and the edge injections for $e$ take the generator of $G_e$ to a generator in $G_w$ and a product of the generators of $G_v$, see Figure~\ref{fig:hc}.
The edge strip is metrized so that closest point projection across the strip is vertical projection in the figure.
Closest point projection from the bottom of the strip to the top of the strip changes distance by a factor of $\frac{1}{\sqrt{2}}$, so the height change across $e$ is $h(v,w)=-\log_2(\frac{1}{\sqrt{2}})=\frac{1}{2}$.

Notice that we get a non-zero height change even though we have identified primitive elements.
\begin{figure}[h]
\labellist
\small
\pinlabel $X_v$ [rt] at 53 719
\pinlabel $X_w$ [rb] at 53 789
\pinlabel $X_e$ [b] at 125 743
\endlabellist
  \centering
  \includegraphics{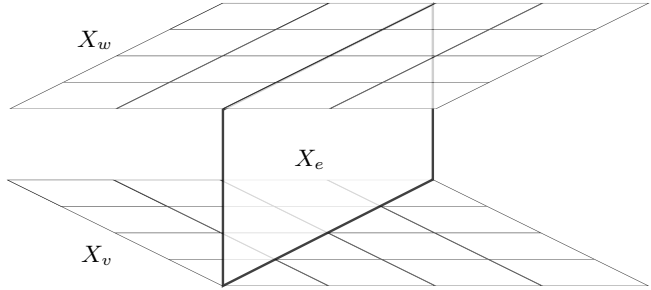}
  \caption{Height change}
  \label{fig:hc}
\end{figure}  
\end{example}

\subsubsection{Geometric Models for Two Tubular Groups}
The coarse Bass-Serre complex, $X'$, will serve as the geometric model for its tubular group.
We will have no further use for the original Bass-Serre complex.

From this point forward, given two tubular groups, $G_i=\pi_1(\Gamma_i)$, for $i=1,2$, $X$ will always refer to the geometric model for $G_1$, and $Y$ will always refer to the geometric model for $G_2$.

\subsubsection{Examples of Tubular Groups}\label{SSS:existence}
We will spend some time discussing a particular family of tubular groups that can be realized as the fundamental group of a graph of groups with one $\mathbb{Z}^2$ vertex group and edge injections contained in three distinct cyclic subgroups. This family is of particular interest for a few reasons:
\begin{itemize}
\item The groups have small enough presentation that they are practical to work with. In \fullref{Ex:bothslopes} we will give a complete quasi-isometry invariant for members of this family.
\item There is only one affine equivalence class of 3-line pattern, so the fact that equivalence classes of edge patterns are quasi-isometry invariants, \fullref{Corollary:patterninvariance}, provides no information.
\item 3-line patterns are rigid, so we will have to worry about preserving height change. This results in a surprisingly delicate quasi-isometry classification.
\item Wise's group and the Brady-Bridson groups all belong to this family. Recall that the Brady-Bridson groups have isoperimetric exponents that form a dense subset of $[2,\infty)$, so we know already that there must be countably many different quasi-isometry classes within the family.
\end{itemize}

\begin{figure}[h]
\labellist
\small \hair 2pt
\pinlabel $v$ [l] at 164 715
\pinlabel $e$ [r] at 99 715
\pinlabel $f$ [l] at 221 715
\pinlabel $\left<a,b\right>$ [r] at 158 715
\pinlabel $a^p$ [bl] at 132 750
\pinlabel $a^q$ [br] at 190 750
\pinlabel $a^rb^s$ [tl] at 131 688
\pinlabel $a^tb^u$ [tr] at 189 688
\endlabellist
\centering
\includegraphics{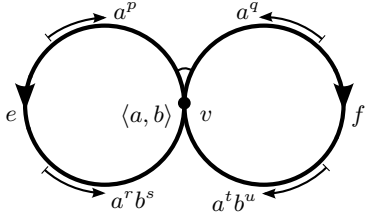}
\caption{A one torus group.}
\label{Fi:OneTorusStandard}
\end{figure}

Let $\Gamma$ be the graph of groups in \fullref{Fi:OneTorusStandard}. The small arc joining edges $e$ and $f$ near their initial endpoints is a visual cue to indicate that the corresponding edge injections map into coarsely equivalent subgroups.
\[G=\pi_1\Gamma =\left<a,b,x,y\mid [a,b]=1,\,x^{-1}a^px=a^rb^s,\, y^{-1}a^qy=a^tb^u\right>\]

Let $F$ be the set of lines through the origin of distinct slopes $0$, $\frac{s}{r}$, and $\frac{u}{t}$.
Let
\[A'=\left(\begin{matrix}
1 & -\frac{1}{2}\frac{ru+st}{su}\\
0 & \frac{\sqrt{3}}{2}\frac{ru-st}{su}
\end{matrix}\right).\]
Let $A=\frac{1}{\sqrt{|\det(A')|}}A'$.
Pull back the metric on $\mathbb{R}^2$ according to the matrix $A$, as in \fullref{P:choosesymmetric}.
This makes the three-line pattern generated by $F$ symmetric.

The height change across $e$ is: \[\lambda=-\log_2 \frac{\left|A(^r_s)\right|}{\left|A(^p_0)\right|}=-\log_2\left|\frac{ru-st}{pu}\right|\]

The height change across $f$ is: \[\mu=-\log_2\frac{\left|A(^{\,t}_u)\right|}{\left|A(^q_0)\right|}=-\log_2\left|\frac{ru-st}{sq}\right|\]

A special case is a Brady-Bridson group $BB(p,r)$ where $0<p=q<r=t$, $s=1$, and $u=-1$.
For such a group the height changes across the two edges are: 
\[\lambda=\mu=-\log_2\frac{2r}{p}\]
The possible values of $\lambda$ are dense in $(-\infty,-1)$.

Brady and Bridson have shown \cite{BraBri00} that the group $BB(p,r)$ has Dehn function: \[n^{-2\lambda}\]
This proved that there are no gaps in the isoperimetric spectrum beyond 2.

\begin{proposition}
Suppose we are given height changes $\lambda$ and $\mu$ and four indices of edge inclusions, $p$, $q$, $m=\gcd(r,s)$, and $n=\gcd(t,u)$.
These can be realized in a tubular group as in \fullref{Fi:OneTorusStandard} if and only if $2^\lambda$ and $2^\mu$ are rational and, as reduced rationals, $\frac{q}{2^\mu n}=\frac{i}{j}$ and $\frac{p}{2^\lambda m}=\frac{i}{k}$ where $i$, $j$, and $k$ are pairwise coprime.
\end{proposition}
\begin{proof}
If these constants did arise from a tubular group, the matrix $A'$ that was used to define the metric on the torus allows us to realize the group $G_v$ as a lattice in $\mathbb{R}^2$.
In this lattice $(1,0)$, $\frac{p}{2^\lambda m}(1/2, \sqrt{3}/2)$, and $\frac{q}{2^\mu n}(-1/2,\sqrt{3}/2)$ must be primitive elements, since these are the images of generators of maximal cyclic subgroups in $G$.
These three elements generate a lattice in which each of them is primitive only if $\frac{p}{2^\lambda m}$ and $\frac{q}{2^\mu n}$ are rational and, as reduced rationals,
 have the same numerators and relatively prime denominators.

For the converse, choose any $g\in \mathbb{Z}$ such that
\[ g\equiv -\frac{1}{j} \text{ mod } k.\]

Choose any $h$ coprime to $g$, $\frac{1+gj}{k}$, and $i$.

Let  $r=mi\left(\frac{1+gj}{k}\right)$, $s=mhj$, $t=ngi$, and $u=nhk$.

With these choices, the graph of groups in \fullref{Fi:OneTorusStandard} gives the desired group.
\end{proof}

\begin{corollary}
For any numbers $\lambda$ and $\mu$ such that $2^\lambda$ and $2^\mu$ are rational,
it is possible to construct a tubular group as in \fullref{Fi:OneTorusStandard} where the height changes across the two edges are $\lambda$ and $\mu$.
\end{corollary}

\subsection{P--sets}\label{S:relationP}
No two depth zero vertex spaces are bounded Hausdorff distance from one another.
However, there is a weaker characterization of closeness of two local spaces according to whether they are close on unbounded subsets.
We define the corresponding relation on the vertices and edges of the Bass-Serre tree:

\begin{definition}
Two elements $x,y\in \mathcal{VE}D\varGamma$ \emph{satisfy Relation P} if $X_x\stackrel{c}{\cap} X_y$ is unbounded.
\end{definition}

Quasi-isometries respect vertex and edge spaces, and preserve boundedness and unboundedness, so Relation P is invariant under quasi-isometry.

\begin{definition}\label{D:Pset}
A \emph{P--set} of $D\varGamma$ is a maximal subset of $\mathcal{VE}D\varGamma$ such that
any two elements satisfy Relation P.
\end{definition} 

\begin{proposition}[Properties of P--sets]\label{P:propertiesofpsets}
\hfill
\begin{enumerate}
\item A P--set is a subtree of $D\varGamma$.\label{Li:psetisconnected}

\item Every depth zero vertex of a P--set is adjacent to infinitely many other vertices in that P--set.\label{Li:infinitevalence}

\item An $n$-line, depth zero vertex belongs to $n$ distinct P--sets.\label{Li:vertinthree}

\item Edges and positive depth vertices belong to exactly one P--set.\label{Li:edgeinone}

\item Any two P--sets are either disjoint or intersect in exactly one vertex, which is necessarily a depth zero vertex.\label{Li:intersectionofpsets}

\end{enumerate}
\end{proposition}
\begin{proof}
Let $e_1$ and $e_2$ be edges in $D\varGamma$ incident to a common depth zero vertex, $v$.

The vertex $v$ satisfies Relation P with either of the $e_i$.

The edges $e_1$ and $e_2$ satisfy Relation P if and only if their edge strips glue to $X_v$ along parallel lines.

The list of properties follows easily from this observation
\end{proof}
The fundamental group $G=\pi_1(\Gamma)$ quasi-acts on $X$.
Quasi-isometries preserve Relation P, so the action of $G$ on $D\varGamma$ induces an action of $G$ on the set of P--sets of $D\varGamma$.

\begin{definition}
Within a P--set $S$, two depth zero vertices $v$ and $w$ are \emph{of the same type} if $\exists g\in \Stab _G(S)$ such that $gv=w$.
\end{definition}

Type is an equivalence relation among depth zero vertices of a fixed P--set.
For an $n$-line depth zero vertex, the equivalence class of the vertex under the action of the whole group splits into at most $n$ vertex types in $S$.

If there is a height change between vertices of the same type in a P--set, the geodesic segment joining them projects to a non-trivial loop in the underlying graph of $\Gamma$.
The group element corresponding to this loop is an infinite order element of the P--set stabilizer.
Iterating the action of this group element gives vertices of the P--set all of the same type occurring at an unbounded set of heights.

If there is zero height change between every pair of vertices of the same type in a P--set then, since there are only finitely many vertex types, the vertices of the P--set occur at only finitely many heights.
Thus, there is bounded height change between any two vertices of the P--set.

For each orbit of P--sets, pick a representative $S_i$ and fix an ordered list of representative $x_{i,j}$ for the vertex types in $S_i$.
Suppose $S$ and $R$ are P--sets, $x$ and $y$ are vertices in $R$, $g\in G$ such that $gR=S$, and $k\in \Stab_G(R)$ such that $kx=y$.
Then $gkg^{-1}\in \Stab_G{S}$ and $gkg^{-1}gx=gkx=gy$, so the action of the group takes vertices of the same type in one P--set to vertices of the same type in another.
Thus, we can fix an ordering of vertex types for each equivalence class of P--set.

We will say that $x\in \mathcal{V}_0S$ is of type $\{[S_i],j\}$ with respect to $S$ if there is some $g\in G$ such that $gS=S_i$ and $gx=x_{i,j}$.

\begin{proposition}\label{P:verttype}
  Suppose $v\in \mathcal{V}_0D\varGamma$, and $R$ and $S$ are P--sets containing $v$.
There is an element $g\in \Stab_G(v)$ such that $gR=S$ if and only if $v$ is of the same type with respect to both $R$ and $S$.
\end{proposition}
\begin{proof}
  Suppose $\exists g\in \Stab_G(v)$ such that $gR=S$.
The P--sets $R$ and $S$ then belong to the same equivalence class; call it $[S_i]$.
Suppose $v\in \{[S_i],j\}$ with respect to $R$ and $v\in\{[S_i],k\}$ with respect to $S$.
Then there are elements $f,h\in G$ with $fR=S_i$, $fv=x_{i,j}$, $hS=S_i$, and $hv=x_{i,k}$.
However, this would mean that $fgh^{-1}\in \Stab_G(S_i)$ and $fgh^{-1}x_{i,k}=x_{i,j}$, but this is only possible if $j=k$.

Conversely, suppose $v\in \{[S_i],j\}$ with respect to both $R$ and $S$.
Then there are elements $f,h\in G$ with $fR=S_i$, $fv=x_{i,j}$, $hS=S_i$, and $hv=x_{i,j}$, so $g=h^{-1}f$ fixes $v$ and takes $R$ to $S$.

\end{proof}

\begin{definition}
The \emph{tree of P--sets,} $T_\varGamma$, of a tubular group, $G=\pi_1(\Gamma)$, is given by:
\begin{itemize}
\item
$\mathcal{V}T_\varGamma=\mathcal{V}_0D\varGamma \amalg \{\text{P--sets of }D\varGamma\}$
\item
Edges are determined by inclusion of vertices in P--sets.
Each edge is assigned length $\frac{1}{2}$.
\end{itemize}
\end{definition}

The collection of P--sets is a cover of $D\varGamma$ by subtrees that intersect pairwise in at most a single vertex. 
The nerve of such a cover is simply connected, and $T_\varGamma$ is a deformation retraction of this nerve. 
So, $T_\varGamma$ is, in fact, a tree.

It will be convenient to extend the notion of height change across an edge to some of the edges of $T_\varGamma$.
If $R$ is a P--set with bounded height change, pick some adjacent depth zero vertex, $v$, of maximal height.
For any vertex $w$ adjacent to $R$, define $h(w,R)=h(w,v)$.

The action of $G$ on $D\varGamma$ induces an action of $G$ on $T_\varGamma$.

\section{Quasi-isometries Between Tubular Groups}
In \fullref{S:allowable} we look at properties of the map of trees of P--sets induced by a quasi-isometry of tubular groups.
A map that satisfies the same properties we call ``allowable.''
In \fullref{ss:allow} we show that an allowable isomorphism of trees of P--sets exists if and only if there is a ``consistent set of strategies,'' and we show that we can decide in finite time whether such a set of strategies exists.
In \fullref{SS:buildingqi} we take a consistent set of strategies and use it to build a quasi-isometry of tubular groups.
Thus, we conclude that the existence of an allowable isomorphism of trees of P--sets is not only a necessary, but also sufficient condition for the existence of a quasi-isometry between tubular groups, and the success or failure of the algorithm for finding a consistent set of strategies determines the existence or non-existence of a quasi-isometry between the tubular groups.
\subsection{Allowable Isomorphisms of Trees of P-Sets}\label{S:allowable}
Suppose $G_1=\pi_1(\Gamma_1)$ and $G_2=\pi_1(\Gamma_2)$ are tubular groups, with trees of P--sets $T_1$ and $T_2$, respectively.

A quasi-isometry, $\psi\from G_1\to G_2$ induces a tree isomorphism $\psi_{\#}\from T_1\to T_2$.

Recall that the P--sets $R_1,\ldots ,R_n$ adjacent to a depth zero vertex $v$ are in one-to-one correspondence with the families of parallel lines forming the edge pattern in the vertex space of $v$. 
To satisfy edge pattern rigidity, the bijection
\[\{R_1,\ldots, R_n\}\to \{\psi_{\#}(R_1),\ldots ,\psi_{\#}(R_n)\}\]
must correspond to an affine equivalence between the edge patterns of $X_v$ and $Y_{\psi_{\#}(v)}$.

If this restriction is satisfied for every depth zero vertex, then $\psi_{\#}$ will be called \emph{locally allowable}.

Define a \emph{rigid component} of $T_i$ to be a connected component of \[T_i\setminus \{2\text{-line, depth zero vertices in }T_i\}\]
An induced isomorphism $\psi_{\#}$ must take 2-line, depth zero vertices to 2-line, depth zero vertices, so $\psi_{\#}$ takes rigid components to rigid components.

\begin{lemma}\label{L:qiimplieshp}
A quasi-isometry between tubular groups induces an isomorphism between their trees of P--sets that is coarsely height preserving on rigid components.
\end{lemma}
\begin{proof}
  \case{Within a P--set Space}{
A P--set space is metrically a warped product of a tree with $\mathbb{R}$, with height as the warping function. 
This also the case for Baumslag-Solitar groups, and the proof that a quasi-isometry is coarsely height preserving within a P--set space is similar to the proof that a quasi-isometry of Baumslag-Solitar groups is coarsely height preserving \cite[Lemma~4.1]{Why01}.

Let $R$ be a P-set in a tubular group $G=\pi_1(\Gamma)$.
In $(X_R, d_{X_R})$, for any vertex, $v\in R$, there is a well defined closest point projection $p_{(R,v)}:X_R\to X_v$.

Suppose $v_0$ and $v_1$ are two vertices in $R$.
In the geodesic of $D\varGamma$ joining $v_0$ and $v_1$, let $e_i$ be the edge incident to $v_i$.
Let $L_{v,e}$ be the line in $X_v$ to which the edge strip for an incident edge $e$ attaches.
Let $x$ and $x'$ be an arbitrary pair of reference points in $L_{v_0,e_0}$.

\[\frac{d_{X_{v_1}}(p(x),p(x'))}{d_{X_{v_0}}(x,x')}=2^{-h(v_0,v_1)}\]

Suppose we have a quasi-isometry, $\phi$, of tubular groups that restricts to a $(K,C)$-quasi-isometry of the $X_{v_i}$.
Let $\phi_\#$ be the induced bijection of depth zero vertices.
Let $S=\phi_\#(R)$ and let $w_i=\phi_\#(v_i)$.
Let $f_i$ be the edge incident to $w_i$ on the tree geodesic joining the $w_i$.
A quasi-isometry coarsely preserves closest point projection. 
That is, the following diagram commutes up to error  $E$ which depends on global constants and on $d_{D\varGamma}(v_0,v_1)$, but not on the $v_i$ or the particular points in $L_{v_0,e_0}$.

\[\begin{CD}
L_{v_0,e_0} @>\phi >>Y_{w_0}\stackrel{c}{\cap}Y_{f_0}@>\stackrel{c}{=}>>L_{w_0,f_0}\\
@Vp_{(R,v_1)} VV @. @VVp_{(S,w_1)}V\\
L_{v_1, e_1}@>\phi>>Y_{w_1}\stackrel{c}{\cap}Y_{f_1}@>\stackrel{c}{=}>>L_{w_1,f_1}
\end{CD}\]

We know that $\phi(x)$ and $\phi(x')$ map to within some uniform distance, $D$, of a line $L\subset Y_{w_0}$ parallel to $L_{w_0,f_0}$.
Thus there is a point $y\in L$ with $d_{Y_{w_0}}(\phi(x),y)\leq D$ and $p_{(S,w_1)}(\phi(x))=p_{(S,w_1)}(y)$.
There is a similar point $y'$ for $\phi(x')$.

The situation is summarized in Figure~\ref{fig:hp}.

\begin{figure}[h]
\labellist
\small
\pinlabel $\phi$ [b] at 284 715
\pinlabel $X_{v_0}$ [br] at 25 631
\pinlabel $X_{v_1}$ [br] at 25 760
\pinlabel $X_{e_1}$ [r] at 90 710 
\pinlabel $X_{e_0}$ [l] at 188 702 
\pinlabel $Y_{w_0}$ [tl] at 530 631
\pinlabel $Y_{w_1}$ [tl] at 530 760
\pinlabel $Y_{f_1}$ [tl] at 476 741 
\pinlabel $Y_{f_0}$ [l] at 503 692

\tiny
\pinlabel $x$ [br] at 117 633
\pinlabel $x'$ [br] at 166 658
\pinlabel $p(x)$ [r] at 107 759
\pinlabel $p(x')$ [r] at 163 787
\pinlabel $L_{v_0,e_0}$ [t] at 93 623
\pinlabel $L_{v_1,e_1}$ [tr] at 93 747
\pinlabel $\phi(x)$ [lt] at 485 626
\pinlabel $y'$ [bl] at 524 663
\pinlabel $y$ [tr] at 455 626
\pinlabel $\phi(x')$ [br] at 505 662
\pinlabel $p\phi(x)$ [br] at 419 754
\pinlabel $p\phi (x')$ [r] at 481 787
\pinlabel $\phi p(x)$ [b] at 419 769
\pinlabel $\phi p(x')$ [l] at 508 791
\pinlabel $L_{w_0,f_0}$ [t] at 408 623
\pinlabel $L_{w_1,f_1}$ [tr] at 407 749
\pinlabel $L$ [lt] at 455 623
\endlabellist  
\centering
  \includegraphics[width=.9\textwidth]{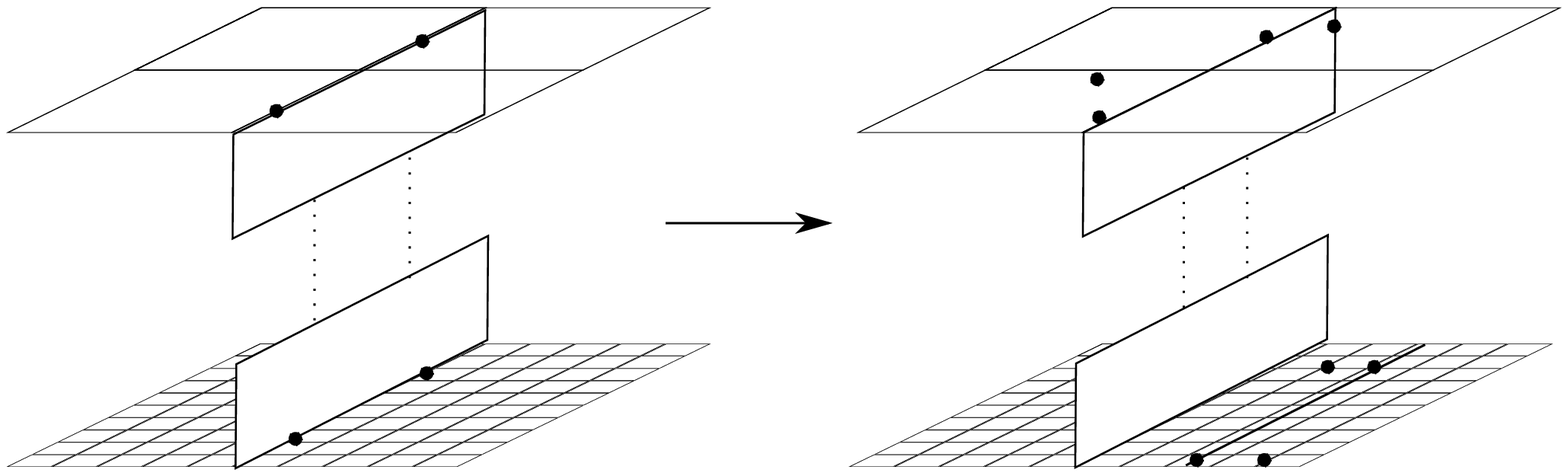}
\caption{}\label{fig:hp}
\end{figure}

\begin{align*}
  2^{-h(w_0,w_1)}&=\frac{d_{Y_{w_1}}\left(p_{(S,w_1)}(y),p_{(S,w_1)}(y')\right)}{d_{Y_{w_0}}\left(y,y'\right)}\\
&\leq\frac{d_{Y_{w_1}}\left(p_{(S,w_1)}(\phi(x)),p_{(S,w_1)}(\phi(x'))\right)}{d_{Y_{w_0}}\left(\phi(x),\phi(x')\right)-2D}\\
&\leq \frac{d_{Y_{w_1}}\left(\phi(p_{(R,v_1)}(x)),\phi(p_{(R,v_1)}(x'))\right)+2E}{d_{Y_{w_0}}\left(\phi(x),\phi(x')\right)-2D}\\
&\leq \frac{Kd_{X_{v_1}}\left(p_{(R,v_1)}(x),p_{(R,v_1)}(x')\right)+C+2E}{\frac{1}{K}d_{X_{v_0}}\left(x,x'\right)-C-2D}
\end{align*}

For any $\delta>0$ we can take $x$ and $x'$ far enough apart so that
\[2^{-h(w_0,w_1)}\leq (1+\delta)K^2\frac{d_{X_{v_1}}\left(p_{(R,v_1)}(x),p_{(R,v_1)}(x')\right)}{d_{X_{v_0}}\left(x,x'\right)}=(1+\delta)K^2\cdot 2^{-h(v_0,v_1)}\]

Thus,
\[h(w_0, w_1)\geq h(v_0,v_1)-2\log_2(K)\]
 The reverse inequality follow from a similar computation, so that height between vertices in a P--set is preserved by quasi-isometries of tubular groups, up to an additive error depending on the quasi-isometry constant.

}

\case{From One P--set Space to Another}{
Let $v_0$, $v_1$, and $v_2$ be depth zero vertices such that $v_i$ and $v_{i+1}$ belong to a common P--set, but $v_0$ and $v_2$ do not.
Assume $v_1$ is a vertex with at least three lines.
For $i=0,1$, let $e_i$ be the edge adjacent to $v_i$ in the tree geodesic joining the $v_i$.
For $j=1,2$, let $e'_j$ be the edge adjacent to $v_i$ in the tree geodesic joining the $v_i$.

Pick reference points $x$ and $x'\in L_{v_0,e_0}$ as in the previous case.
Projecting these points to $X_{v_2}$ would not give us much information because $X_{v_0}\stackrel{c}{\cap}X_{v_2}$ is a point in $X_{v_1}$.
We can not do two separate height change calculations within P--sets because we need the height error bound to be independent of the number of P--sets traversed.
Instead we link these calculations together by exploiting the fact that the restriction of a quasi-isometry to a vertex space with at least three lines expands by the same factor in all directions.

Pick reference points $z$ and $z'\in L_{v_1,e'_1}$ such that 
\[d_{X_{v_1}}(z,z')=d_{X_{v_1}}(p_{(R,v_1)}(x), p_{(R,v_1)}(x'))\]

Since $\phi|_{X_{v_1}}$ stretches by the same amount in every direction, we have
 \[ d_{Y_{w_1}}\left(\phi(z),\phi(z')\right)\leq d_{Y_{w_1}}\left(\phi(p_{(R,v_1)}(x)),\phi(p_{(R,v_1)}(x'))\right)+D_1\]
Thus, for any $\delta_1>0$ we can take $x$ and $x'$ far enough apart so that 
 \[ d_{Y_{w_1}}\left(\phi(z),\phi(z')\right)\leq (1+\delta_1)\cdot d_{Y_{w_1}}\left(\phi(p_{(R,v_1)}(x)),\phi(p_{(R,v_1)}(x'))\right)\]

Now we can chain together the height change calculations from the different P--sets and get cancellation at the intermediate vertices. 
We find that for any $\delta_2>0$, if $x$ and $x'$ are sufficiently far apart then 
\[2^{-h(w_0,w_2)}\leq (1+\delta_2)K^2\cdot 2^{-h(v_0,v_2)}\]

So as before we find that \[h(w_0, w_2)\geq h(v_0,v_2)-2\log_2(K)\]
A similar computation produces the reverse inequality.
}
\end{proof}

\begin{definition}\label{Def:allowable}
A tree isomorphism $\phi\from T_1\to T_2$ is \emph{allowable} if it is locally allowable and coarsely height preserving on rigid components.
\end{definition}

\fullref{L:qiimplieshp} and the preceding material then give the following corollary.
\begin{corollary}\label{Corollary:qiimpliesallowable}
  Quasi-isometries of tubular groups induce allowable isomorphisms on their trees of P--sets.
\end{corollary}

\subsection{Finding Allowable Isomorphisms}\label{ss:allow}
Let $[R_1],\ldots ,[R_m]$ be the equivalence classes of P--sets in $D\varGamma_1$ under the action by $G_1$.

Let $[S_1],\ldots ,[S_n]$ be the equivalence classes of P--sets in $D\varGamma_2$ under the action by $G_2$.

Let $\rho(a)$ be the number of vertex types in $[R_a]$ and let $\sigma(b)$ be the number of vertex types in $[S_b]$.

A \emph{match} is a pair $([R],[S])$ consisting of an equivalence class of P--set from $D\varGamma_1$ and one from $D\varGamma_2$.

An \emph{extension} for $([R_a],[S_b])$ is a $\rho(a)\times \sigma(b)$ matrix, $(e_{ij})$, with entries described below.
An extension must have at least one non-zero entry in each row and in each column.

Recall that we have chosen representatives $x_{a,i}\in \{[R_a],i\}$ with respect to $R_a$ and $y_{b,j}\in \{[S_b],j\}$ with respect to $S_b$.
Each $x_{a,i}$ belongs to some P--sets $R^{a,i,l}$ (one of which is $R_a$) where $l=1\dots \text{number of lines in $x_{a,i}$}$.
Consider the set, possibly with repeated entries, $\{\{[R^{a,i,l}],k_l\}\}_{l}$, where $x_{a,i}$ is of type $\{[R^{a,i,l}],k_l\}$ with respect to $R^{a,i,l}$. This is the set of types that $x_{a,i}$ takes with respect to the P--sets containing it. 

Similarly, consider the set $\{\{[S^{b,j,l}],k_l'\}\}_l$, where $y_{b,j}$ is of type $\{[S^{b,j,l}],k_l'\}$ with respect to $S^{b,j,l}$.

A linear equivalence between the edge patterns in $X_{x_{a,i}}$ and $Y_{y_{b,j}}$ gives a bijection between P--sets adjacent to $x_{a,i}$ and $y_{b,j}$, which in turn gives a bijection of these sets of types.
We will let $e_{ij}$ be any bijection of these sets that can be induced by a linear equivalence of edge patterns and that includes $\{[R_a],i\}\to \{[S_b],j\}$.
\begin{remark}
If we knew that for any parallel family of edge lines in a vertex space, the subgroup of the vertex stabilizer that fixes that family also fixes each of the other families, then the set of types that that depth zero vertex takes has no duplicate entries. 
This would be true, for instance, if the tubular group were torsion free.
Specifying a type picks out a parallel family of lines in the vertex space, and we could put an equivariant numbering on these families.
Then, instead of a bijection of sets of types as we defined above, we could just associate to a linear equivalence of edge patterns a permutation describing which families go to which families.
This may fail for groups with torsion, but \fullref{P:verttype} shows that it is just the ordering of the sets that is ambiguous.
In either case, there are at most two distinct bijections to consider.
\end{remark}
If there are no such bijections then $e_{ij}=0$.

Notice that the number of possible extensions for $([R_a],[S_b])$ is bounded above by 
\[3^{\rho(a) \sigma(b)}.\]

An extension provides instructions for extending a tree isomorphism in a locally allowable way.
Suppose $R\in T_1$ and $S\in T_2$ are P--sets and we have constructed a map sending $R$ to $S$.
To extend the map, choose any bijection $\mathcal{V}_0R\to \mathcal{V}_0S$ such that vertices of type $\{[R],i\}$ map to vertices of type $\{[S],j\}$ if and only if $e_{ij}\neq 0$.
The value of $e_{ij}$ then tells how to extend to the next level of P--sets.
For the purposes of tree isomorphisms, the particular bijection  $\mathcal{V}_0R\to \mathcal{V}_0S$ does not matter, but more care will be required when building quasi-isometries in \fullref{SS:buildingqi}.

Conversely, given an allowable isomorphism $\phi \from N_1(R)\to N_1(S)$, it is possible to ``read off'' an extension for $([R],[S])$.
Set $e_{ij}=0$ if $\phi$ does not take any vertices of type $\{[R],i\}$ to vertices of type $\{[S],j\}$.
Otherwise, choose some representative vertex $v_{i,j}$ with $v_{i,j}\in \{[R],i\}$ and $\phi(v_{i,j})\in \{[S],j\}$.
Set $e_{ij}$ equal to the bijection of vertex types of $v_{i,j}$ to vertex types of $\phi(v_{i,j})$ induced by $\phi$.
Then the matrix $(e_{ij})$ is an extension for $([R],[S])$.

\subsubsection{Strategies}
A \emph{strategy,} $\mathcal{S}$, for $([R],[S])$ consists of:
\begin{enumerate}
\item a root vertex, $\rm{ROOT}(\mathcal{S})$, which is labeled by $([R],[S])$,
\item an extension $\mathcal{E}=(e_{ij})$ for $([R],[S])$,
\item a collection of terminal vertices corresponding to the set $\{(\mathcal{M}, E)\}$ of induced matches with height errors coming from $\mathcal{E}$.
\end{enumerate}

The \emph{label} of a terminal vertex is the match associated to it.
A strategy can be simplified by considering at most three terminal vertices for each label: one with maximum height error, one with minimum height error, and one with undefined height error.
Thus, we can assume that the number of terminal vertices of any strategy is at most $3mn$.

A strategy records how the boundary of a neighborhood of a P--set $R$ in $T_1$ maps to the boundary of a neighborhood of $S$ in $T_2$ when mapped according to an extension $\mathcal{E}$ for $([R],[S])$.

If $\mathcal{S}$ is a strategy, define \[\rm{TERM}(\mathcal{S})=\text{$\{$terminal vertices of $\mathcal{S}\}$}\]

A \emph{set of strategies} will consist of a graph, each of whose vertices is labeled by a match and two (not necessarily distinct) strategies for the match, one called the positive strategy and the other called the negative strategy.
There will be at most one vertex labeled by a given match, so at most $mn$ vertices.

For each strategy in the set, add an edge to the graph from the root vertex of the strategy to each terminal vertex of the strategy.
Label the edge by the sign of the strategy and the height error, if defined, for the appropriate terminal vertex.

\subsubsection{Consistency}
We need to check that building according to the set of strategies will not create unbounded height error.
It is necessary that each match consists of two P--sets of bounded height change or two P--sets of unbounded height change, but we also must control accumulation of height error created by the strategies.

For P--sets of bounded height error, the type of a vertex determines its relative height.
Choosing an extension therefore determines height errors in a neighborhood of the P--set.

We do not need to worry about height change when it comes to P--sets of unbounded height change.
We will see in \fullref{L:qiofUnboundedPsets} that the presence of depth zero vertices at an unbounded set of heights allows us to correct any height error.

Let $\mathcal{M}_1,\ldots $ be a list of the matches of P--sets of bounded height change that label vertices of $\mathcal{G}$.
For each $i$, add to the set of inequalities:
\[L_i\leq E_i \leq U_i\]
Let $\mathcal{S}_i^+$ and $\mathcal{S}_i^-$ be the positive and negative strategies chosen for $\mathcal{M}_i$.

Suppose $v\in \rm{TERM}(\mathcal{S}_i^+)$ has a defined height error, and the label of $v$ is $\mathcal{M}_j$.
Let $E$ be the height error of $v$ in $\mathcal{S}_i^+$.
Add the following inequalities to the set:
\begin{gather*}
M_i + E \geq L_j\\
U_i + E \leq U_j
\end{gather*}
Suppose $v\in \rm{TERM}(\mathcal{S}_i^-)$ has a defined height error, and the label of $v$ is $\mathcal{M}_j$.
Let $E$ be the height error of $v$ in $\mathcal{S}_I^-$.
Add the following inequalities to the set:
\begin{gather*}
M_i + E \leq U_j\\
L_i + E \geq L_j
\end{gather*}

If the system of inequalities has a solution then height error is well controlled.
The $L_i$ and $U_i$ provide bounds.

If such a system of inequalities has a solution, then it has a solution such that all the $U_i$ are positive, all the $L_i$ are negative, and all the $M_i$ are zero.
Furthermore, given any fixed $B$, there is a solution such that, for all $i$, $U_i>B$ and $L_i<-B$.

Note that if a positive strategy for a match $\mathcal{M}_i$ adds an edge with negative height error that leads back to $\mathcal{M}_i$, then the set of strategies can not be consistent.
Similarly, a negative strategy should not create a length one loop with a positive height error.

The next lemma shows that from an allowable isomorphism we can derive a consistent set of strategies. 
The idea of the proof is to consider the places where the isomorphism has the worst height errors and choose strategies according to what happens in neighborhoods of those bad places.

\begin{lemma}[Deriving a Set of Strategies]\label{L:Success}
If there is an allowable isomorphism between $T_1$ and $T_2$, then a consistent set of strategies exists.
\end{lemma}
\begin{proof}
Suppose $\phi$ is an allowable isomorphism $T_{1}\to T_{2}$ of trees of P--sets.

For each match $([R],[\phi(R)])$ with $[R]$ and $[\phi(R)]$ of unbounded height change, pick some representative $R$ and read off a strategy.
This strategy suffices for both the positive and negative strategy of $([R],[\phi(R)])$.

If there is a P--set $R$ of bounded height change such that every depth zero vertex of $R$ has two lines then read off a strategy for $R$.
This strategy suffices for both the positive and negative strategy of $([R],[\phi(R)])$.

The only matches left are pairs of P--sets of bounded height change that contain at least one depth zero vertex with at least 3 lines.
Pick such a P--set, $R$, and depth zero vertex, $v_0$.

Every rigid component, $\mathcal{C}$, of $T_1$ contains a vertex, $v_{\mathcal{C}}$, closest to $v_0$.
If $v_{\mathcal{C}}\neq v_0$ then $v_{\mathcal{C}}$ is necessarily a 2-line, depth zero vertex.

For any depth zero vertex, $w$, in the rigid component $\mathcal{C}$, define $h(w)=h(v_\mathcal{C},w)$.

Make the corresponding definitions for heights of depth zero vertices of $T_2$.

Let $\mathcal{M}_1$, \ldots, $\mathcal{M}_a$ be a list of matches that occur as $([R], [\phi(R)])$, with $R$ and $\phi(R)$ of bounded height change and having a depth zero vertex with at least 3 lines.
Let 
\begin{gather*}
U_i=\sup\{h(R)-h(\phi(R))\mid ([R],[\phi(R)])=\mathcal{M}_i\}\\
L_i=\inf\{h(R)-h(\phi(R))\mid ([R],[\phi(R)])=\mathcal{M}_i\}
\end{gather*}
These quantities exist since $\phi$ is uniformly coarsely height preserving on rigid components.

For any match there are only finitely many possible strategies.
For each $1\leq i \leq a$, if $U_i$ is achieved, pick an $R$ with $U_i=h(R)-h(\phi(R))$ and  $([R],[\phi(R)])=\mathcal{M}_i$ and read off a strategy $\mathcal{S}_i^-$ for $([R],[\phi(R)])$.
If $U_i$ is not achieved, let $\mathcal{S}_i^-$ be any strategy that occurs for $R$ with height error arbitrarily close to $U_i$.
These will be the negative strategies.

Pick positive strategies $\mathcal{S}_i^+$ analogously, using the $L_i$.

The set of strategies chosen in this way is consistent, although the height error bounds may be larger than for $\phi$.
\end{proof}

\begin{theorem}\label{T:Algorithm}
There is an algorithm that in finite time either produces a consistent set of strategies for two trees of P--sets or decides that no such set exists.
Furthermore, the algorithm is guaranteed to succeed if the trees of P--sets came from quasi-isometric tubular groups.
\end{theorem}
\begin{proof}
A set of strategies is completely determined by the matches in it and the choices of extensions for those matches.
If $T_1$ and $T_2$ are as above, then there are at most $mn$ matches and each extension has at most $\max \{\rho(a)\} \times \max \{\sigma(b)\}$ entries.
Each entry takes one of at most 3 values, so we have at most
\[3^{2mn\times \max \{\rho(a)\} \times \max \{\sigma(b)\}}\]
possible candidate sets of strategies.

Enumerate this list of candidates and check if there is a candidate that is actually a consistent set of strategies.

In light of \fullref{Corollary:qiimpliesallowable} and \fullref{L:Success}, success is guaranteed if the trees came from quasi-isometric groups.
\end{proof}

\subsection{Building Quasi-isometries}\label{SS:buildingqi}
In this subsection we show how to build a quasi-isometry of tubular groups from a consistent set of strategies.
In \fullref{ss:trees} we prove some auxiliary lemmas about quasi-isometries of trees. 
In \fullref{ss:pieces} we build quasi-isometries of P--set spaces.
In \fullref{ss:equiv} we assemble quasi-isometries of P--set spaces to get a quasi-isometry of tubular groups.
\subsubsection{Some Lemmas on Trees}\label{ss:trees}
We leave the world of tubular groups for a moment to prove some lemmas about quasi-isometries of trees that we will need to build quasi-isometries of P--set spaces.

Let $\FIN(V)$ denote the set of all finite subsets of $V$.

\begin{refthm}[Hall's Selection Theorem]
Let $\omega:V \to \FIN(W)$.
There exists $\phi:V\to W$ an injection with $\phi(v)\in \omega(v)$ if and only if for all $S\in\FIN(V)$, $|S|\leq |\cup_{s\in S}\omega(s)|$
\end{refthm}

A tree will always mean a simplicial tree with edges of length 1.
A line $r$ in a tree $V$ means an isometric embedding $r\from \mathbb{R}\into V$ that takes integers to vertices.

A tree is $\delta$-\emph{bushy} if every vertex is within distance $\delta$ of a vertex whose complement has at least three unbounded components.

A line $r$ has \emph{$C$-coarse slope} $m$ if, for all $a$, $b\in\mathbb{R}$, 
\[\left| h(r(a),r(b))-m|b-a|\right|\leq C\]

Let $E_r(a,b)$ denote the number of edges incident to the segment of $r$ from $r(a)$ to $r(b)$.
A line has \emph{$C$-coarse edge density} $\rho$ if, for all $a$, $b\in\mathbb{R}$, 
\[\left| E_r(a,b)-\rho |b-a|\right|\leq C\]

A \emph{lamination} of $V$ is a family of lines such that every vertex belongs to exactly one of the lines.

Whyte proved \cite{Why01} the quasi-isometry classification of graphs of $\mathbb{Z}$'s by matching laminations of their Bass-Serre trees.
In \fullref{L:generalizedBS} we will give a generalization of Whyte's argument, but first we will give an outline of his proof.
The goal is to build a coarsely height preserving quasi-isometry between two trees, $V_1$ and $V_2$.

\begin{sketch}
\hfill
\begin{steps}
\step
Reduce to the case of homogeneous trees.

\step
For any sufficiently small $\beta$ there is a lamination of $V_i$ by lines of coarse slope $\beta$.
Since the $V_i$ are homogeneous, any line in the tree has coarse edge density $\rho_i$ equal to the valence minus two.
Choose sufficiently small $\beta_i$ so that: 
\[\frac{\beta_1}{\beta_2}=\frac{\rho_1}{\rho_2}\]
Choose laminations of the $V_i$ by lines of slope $\beta_i$.

\step
The quasi-isometry will be built line-by-line. Given a line $r_i$ from each lamination, there is an obvious coarsely height preserving quasi-isometry: rescale the line from $V_1$ by a factor of $\frac{\beta_1}{\beta_2}$.

\step
Since the ratio of the slopes was equal to the ratio of the edge densities, a segment of $r_1$ has approximately the same number of incident edges as its image in $r_2$.
This allows us to produce a height respecting  matching of lines of the lamination of $V_1$ adjacent to $r_1$ to lines of the lamination of $V_2$ adjacent to $r_2$.

\step
Induct.
\end{steps}
\end{sketch}

We wish to generalize this argument to trees that are not Bass-Serre trees of a group.
Let $\Gamma$ be a finite, connected graph with directed edges.
The graph $\Gamma$ may have edges with the same initial and terminal vertex, and may have multiple edges between a pair of vertices.
Associate a height change to each edge.
Suppose there is a loop in $\Gamma$ such that the sum of the height changes across edges of the loop is strictly positive.

Consider a bounded valence tree with directed edges covering $\Gamma$, $V\stackrel{p}{\to} \Gamma$, such that for every vertex $v\in V$, the edges coming in to $v$ cover the edges coming in to $p(v)$ at least two to one, and similarly for the outgoing edges.
If $\Gamma$ is not just a single vertex with a single edge, then this condition can be relaxed slightly.
If there is an edge $e$ and a vertex $v$ in $\Gamma$ such that $v$ is the initial and terminal vertex for $e$, then at a vertex in $p^{-1}(v)$ there need be only one incoming and one outgoing edge covering $e$.

The height change of an edge in $V$ is the height change of its image in $\Gamma$.

To prove there is a coarsely height preserving quasi-isometry between two such trees we follow Whyte's outline. 
We can not reduce to homogeneous trees, so there is no reason to believe that we can choose a lamination whose lines have a well defined edge density.
However, choosing lines with well defined slopes and edge densities really amounts to saying that there are uniform proportionality constants such that for any segment of a line of the lamination, the length of the segment is proportional to the height change along the segment and to the number of edges incident to the segment.
We do not need this much to make Whyte's argument work. Really we need two conditions. We need bounds on the ratio of height change to length of a segment so that we can build height preserving quasi-isometries along the lines, and we need the ratio of height change to number of incident edges to be uniform on all segments of lines of the lamination so that we can match adjacent lines.

\begin{lemma}\label{L:generalizedBS}
Let $V\stackrel{p}{\to}\Gamma$ and $W\stackrel{q}{\to}\Gamma'$ be two trees as described above.
There is a coarsely height preserving quasi-isometry between $V$ and $W$.
Furthermore, the quasi-isometry and coarseness constants can be bounded in terms of information from $\Gamma$ and $\Gamma'$ and the valence bounds for $V$ and $W$.
\end{lemma}
\begin{proof}

Assume $W$ is the $(2,2)$-homogeneous tree, that is, the tree that has at every vertex two edges that increase height by one, and two edges that decrease height by one.
For the general statement of the Lemma it suffices to compose two instances of this special case.

Let $M$ be the maximum height change across an edge of $\Gamma$.

Let $\delta$ be the diameter of $\Gamma$.

Let $N$ be the number of vertices of $\Gamma$.

Let $L$ be the maximum valence of a vertex in $V$.

\begin{steps}
\step
Take a maximal subtree of $\Gamma$ with a basepoint and a family of lifts of this subtree to $V$ such that every vertex of $V$ belongs to exactly one lift in the family.
The covering is locally 2 to 1 on edges of the maximal subtree, so such a family exists.
Collapsing these subtrees is a $(\delta+1, \delta)$-quasi-isometry to a tree, $V'$.
We can make it $2\delta M$-coarsely height preserving by declaring the height change between two vertices of $V'$ to be the height change between lifts of the basepoint in the two preimages in $V$.
The maximum height change across an edge of $V'$ is at most $M'=(2\delta+1)M$.
By assumption, there was a loop in $\Gamma$ that strictly increased height, so every vertex of the tree has at least one edge that increases height and one that decreases height.
Furthermore, there is no vertex that has exactly one edge that decreases height and all the rest increase height, or vice, versa.
A vertex with only one edge that decreases height or only one edge that increases height must also have edges with zero height change.
If any such vertices exist, then further collapse a collection of disjoint, zero-height-change edges so that every vertex in the tree has at least two incident edges that increase height and at least two that decrease height.

Vertices of $V'$ have valence between 4 and $2NL$.

\step
Let $\rho=2NL-2\geq 2$.
Let $0<\alpha \leq 1$ be some number such that every vertex of $V'$ has at least two edges that increase height by at least $\alpha$ and two edges that decrease height by at least $\alpha$.
Such an $\alpha$ exists since there were only finitely many possible height changes across an edge of $V'$.

Let $\beta=\frac{2\alpha}{\rho}\leq 1$.
There is a lamination of $W$ by lines of 2-coarse slope $\beta$. 
Such a lamination can be built inductively. Start at some vertex and choose an incident edge to start a line. To continue a line choose either an edge that increases height or decreases height as appropriate to keep the slope of the line as close as possible to $\beta$. 
An appropriate edge is always available since at every vertex there are two edges that increase height and two edges that decrease height.

Repeat this process for every vertex adjacent to the lines that have already been built, etc. 

\begin{claim}
There exists a constant $J>0$ and a lamination of $V'$ by lines $r$ such that, for all $m$, $n\in\mathbb{R}$,
\[\left| h(r(m),r(n))-\frac{\alpha}{\rho}E_r(m,n)\right| \leq J\]
\end{claim}

\step
Assuming the claim, we build a quasi-isometry line-by-line using the laminations of $V'$ and $W$.

For a line $r$ in the lamination of $V'$, \[2(|b-a|+1)\leq E_r(a,b) \leq \rho(|b-a|+1)\]
so \[\beta |b-a|+\beta-J\leq h(r(a),r(b)) \leq \alpha |b-a| +\alpha +J\]

It is sufficient to define a quasi-isometry on the integer points.

Assume that $\phi(r(0))= r'(0)$. 
Let these points be the basepoints of their respective trees, and define height in each tree relative to the basepoint, that is, $h(r(a))=h(r(0),r(a))$.

For $z\in\mathbb{Z}$ pick a $z'\in\mathbb{Z}$ so that \[|h(r'(z'))-h(r(z))|\leq \frac{1}{2}\] and define $\phi(r(z))=r'(z')$.
Such a $z'$ exists because $r'$ has nonzero slope and height changes by $\pm 1$ across each edge of $W$.

For any $z_1$, $z_2\in \mathbb{Z}$,
\[h(r(z_1))-h(r(z_2))-1\leq h(r'(z_1'))-h(r'(z_2'))\leq h(r(z_1))-h(r(z_2))+1\]

Combining these inequalities with the bounds on height change in terms of segment length for $r$ and the fact that $r'$ has 2-coarse slope $\beta$, we find:
\[\beta|z_1-z_2|+\beta-J-1\leq \beta |z_1'-z_2'|+2\]
and
\[\beta |z_1'-z_2'|-2\leq \alpha |z_1-z_2|+\alpha +J+1\]
Thus, $\phi$ is a coarsely height preserving quasi-isometry.

\step

For $x<y\in\mathbb{R}$, let $\mathcal{V}_1(r)[x,y]$ be the set of vertices adjacent to $r$ with height between $x$ and $y$.

We need a $K>0$ so that for all $x$ and $y$, \[\left|\mathcal{V}_1(r)[x,y]\right|\leq \left|\mathcal{V}_1(r')[x-K,y+K]\right|\] and vice versa.

Then Hall's Selection Theorem gives injections between vertices adjacent to $r$ and vertices adjacent to $r'$ that are $K$-coarsely height preserving, and the Schroeder-Bernstein Theorem applied to these injections gives a bijection.

Suppose $a<b$ are such that $h(r(a))=x$ and $h(r(b))=y$. Then:
\[\beta |b-a|+\beta-J\leq |y-x| \leq \alpha |b-a| +\alpha +J\]

Let $c=\frac{M'+J}{\beta}$. Points on $r$ distance at least $c$ apart differ in height by more than $M'$,
so if $z<a-c$ then $x-h(r(z))=h(r(a))-h(r(z))>M'$. Since edges of $V'$ change height by at most $M'$, every vertex adjacent to $r(z)$ has height strictly less than $x$, so no vertices adjacent to $r(z)$ belong to $\mathcal{V}_1(r)[x,y]$.

By similar arguments, no vertices adjacent to $r(z)$ for $z>b+c$ are in $\mathcal{V}_1(r)[x,y]$, and all vertices adjacent to $r(z)$ for $z\in [a+c,b-c]$ are in $\mathcal{V}_1(r)[x,y]$.

Therefore:
\begin{align*}
E_r(a+c,b-c)&\leq \left|\mathcal{V}_1(r)[x,y]\right|\leq E_r(a-c,b+c) \\
&<E_r(a-c,a+c)+E_r(a+c,b-c)+E_r(b-c,b+c)
\end{align*}

We can estimate these quantities:
\[\left|\frac{\beta}{2}E_r(a+c,b-c)-(y-x)\right|\leq 2\alpha(c+1)+3J\]

\[\frac{\beta}{2}E_r(a-c,a+c)\leq \alpha(2c+1)+2J\]

\[\frac{\beta}{2}E_r(b-c,b+c)\leq \alpha(2c+1)+2J\]

This gives us:
\[\left|\frac{2}{\beta}(y-x)- \left|\mathcal{V}_1(r)[x,y]\right|\right|\leq \frac{4}{\beta}\left(\alpha(3c+2)+\frac{7}{2}J\right)\]

A similar computation for $r'$ shows:
\[\left|\frac{2}{\beta}(y-x)- \left|\mathcal{V}_1(r')[x,y]\right|\right|\leq \frac{20}{\beta} \]

Let $A=\max\{\frac{4}{\beta}\left(\alpha(3c+2)+\frac{7}{2}J\right),\frac{20}{\beta}\}$, and let $K=\frac{A\beta}{2}$.

Then 
\begin{align*}
  \left|\mathcal{V}_1(r')[x,y]\right|&\leq \frac{2}{\beta}(y-x)+A\\
&= \frac{2}{\beta}(y-x+2K)-A\\
&\leq \left|\mathcal{V}_1(r)[x-K,y+K]\right|
\end{align*}

Similarly: \[  \left|\mathcal{V}_1(r)[x,y]\right|\leq \left|\mathcal{V}_1(r')[x-K,y+K]\right|\]

\step Induct.

\end{steps}
\bigskip
\paragraph*{Proof of the Claim}
It is sufficient to take $J=2J'$ where $J'=M'+\alpha$.

Let $r(0)$ be some vertex.
There is an edge incident to $r(0)$ that increases height by at least $\alpha$.
Follow this edge to $r(1)$.

\[4\leq E_r(0,1)\leq 2\rho\]
\[\alpha\leq h(r(0),r(1)) \leq M'\]
So:

\[-J'<\alpha-2\rho\frac{\beta}{2} \leq h(r(0),r(1))-\frac{\beta}{2}E_r(0,1)\leq M'-4\frac{\beta}{2}<J'\]

\step
Continue by induction.
Suppose that $r$ has been extended to $r(n)$.
If \[0< h(r(0),r(n))-\frac{\beta}{2}E_r(0,n)\leq J'\]
then follow an edge that decreases height by at least $\alpha$ to get to $r(n+1)$.
There are at least two edges incident to $r(n)$ that decrease height by at least $\alpha$, so it is possible to extend without backtracking.
\begin{align*}
-J'&< h(r(0),r(n))-\frac{\beta}{2}E_r(0,n)-M'-\alpha \\
&\leq h(r(0),r(n+1))-\frac{\beta}{2}E_r(0,n+1)\\
&\leq h(r(0),r(n))-\frac{\beta}{2}E_r(0,n)-\alpha-\beta<J'
\end{align*}

Conversely, suppose \[-J'\leq h(r(0),r(n))-\frac{\beta}{2}E_r(0,n)\leq 0 \]
Without backtracking, follow an edge that increases height by at least $\alpha$ to get to $r(n+1)$.
\begin{align*}
  -J'&\leq h(r(0),r(n)) -\frac{\beta}{2}E_r(0,n)+\alpha-\alpha\\
&\leq h(r(0),r(n+1))-\frac{\beta}{2}E_r(0,n+1)\\
&\leq h(r(0),r(n))-\frac{\beta}{2}E_r(0,n)+M'-\beta<J'
\end{align*}

Repeat this process to build a ray from $r(0)$ such that 
\[\left| h(r(0),r(-n))+\frac{\beta}{2}E_r(0,-n)\right| \leq J'\]

The union of the two rays gives the desired line.

Take some vertex adjacent to this line.
This vertex has at least two edges that increase slope, and two edges that increase slope.
Only one of these edges is the one that connects the vertex to the line, so we can repeat the construction to get another line through this vertex disjoint from the first line.
We can continue to extend in this way to get a lamination of $V'$.
\end{proof}

It will be important that the quasi-isometries we build are bijections on the vertex sets of the trees.
Since the trees we are working with are non-amenable this can always be arranged.
The $n=1$ case of the following lemma is a special case of Theorem~4.1 of \cite{Why99}.
\begin{lemma}\label{L:3valent}
Let $V$ and $W$ be $\delta$-bushy trees, both with variable valence at most $b$.
Suppose $\psi:V\to W$ is a quasi-isometry.
Given finitely many disjoint, $\delta$-coarsely dense subsets $V_1,\ldots ,V_n\subset \mathcal{V}V$ and $W_1,\ldots, W_n\subset \mathcal{V}W$,
there exists $\chi:\cup V_i\to \cup W_i$ a bilipschitz bijection that respects the partitions, $\chi(V_i)\subset W_i$ for all $i$.
The map $\chi$ extends to a quasi-isometry $V\to W$.
Furthermore, $\chi$ is bounded distance from $\psi$, and the quasi-isometry constants depend only on $b$, $\delta$, and the quasi-isometry constants of $\psi$.
\end{lemma}
\begin{proof}
The trees $V$ and $W$ are both quasi-isometric to the Cayley graph of the free group of rank 2, so they are both non-amenable.
Suppose $\psi$ is a $(\lambda, \epsilon)$-quasi-isometry. 

Since each of the $V_i$ and $W_i$ are $\delta$-coarsely dense, $\psi$ can be changed by distance at most $\delta$ so that for each $i$ we have a $(\lambda, 2\delta (\lambda +\epsilon))$-quasi-isometry $\psi_i=\psi |_{V_i}\from V_i\to W_i$.
The $V_i$ and $W_i$ are also non-amenable, since they are coarsely dense subsets of non-amenable sets.

For each $\psi_i$ apply Theorem~4.1 of \cite{Why99} to get a bilipschitz bijection $\chi_i\from V_i\to W_i$.
The $V_i$ are disjoint, so we can assemble the $\chi_i$ to get a bilipschitz bijection $\chi\from \cup V_i \to \cup W_i$ that respects the partitions.
Finally, since $\cup V_i$ is coarsely dense in $V$, $\chi$ can be extended to a give a quasi-isometry $V\to W$. 

\end{proof}

\subsubsection{The Pieces}\label{ss:pieces}
A \emph{P--set space} is the preimage in the model space of a P--set in the Bass-Serre tree.
In \fullref{L:qiofPsets} and \fullref{L:qiofUnboundedPsets} we construct quasi-isometries of P--set spaces.
In \fullref{T:Construction} these are pieced together to give a quasi-isometry of tubular groups.

There are two complications to consider when building quasi-isometries between P--set spaces.
First, depth zero vertex spaces belong to several different P--set spaces.
We must take care to build quasi-isometries of P--set spaces that can later be pieced together into a quasi-isometry of tubular groups.
Second, to make use of \fullref{T:Algorithm} we need to build according to a given extension, which puts type restrictions on which vertex spaces map to each other.

Let $R$ be a P--set in $D\varGamma$.
Let $v_0$ be some depth zero vertex of $R$.
Height for other points of $R$ can then be defined relative to $v_0$.

Put standard Euclidean coordinates on $X_{v_0}$ in such a way that the edges of $R$ glue onto $X_{v_0}$ along lines of the form ${t}\times \mathbb{R}$.

The second coordinate can be projected over all of $X_R$, giving coordinates on $X_R$ of the form $\text{tree} \times \mathbb{R}$, with the $\mathbb{R}$ factor scaled according $2^{-\text{height}}$.
With these coordinates, an edge space or a positive depth vertex space is just $(t,\mathbb{R})$, for some $t$, and the length
of $(t,[a,b])$ is $2^{-h(t)}|b-a|$.

If $v$ is a depth zero vertex of $R$, parameterize the direction orthogonal to the lines of edge attachments so that the edge leading back to $v_0$ is at coordinate 0 and so that $d_{X_{v}}=2^{-h(v)}d_{\mathbb{R}\times \mathbb{R}}$.

Define a map $B\from R\to \text{tree}$ by lifting to $X_R$ and projecting to zero in the second coordinate.
In $B(R)$ the depth zero vertices are ``blown up'' into lines.
Edges of $B(R)$ contained in a depth zero vertex space will be called \emph{horizontal}.
Conversely, edges of $B(R)$ that cross an edge space will be called \emph{vertical}.

A \emph{horizontal component} of $B(R)$ is the blow up of some depth zero vertex.
A \emph{vertical component} is a connected component of the complement of the interiors of the horizontal edges.

We will build a quasi-isometry of $B(R)$ that can be extended to a quasi-isometry of the P--set spaces.
At first glance the situation here is very similar to the situation described in the previous section, we want a height preserving quasi-isometry of trees and we (almost) have a lamination of the trees by the collection of horizontal components.
However, there are some important differences, and the techniques of \fullref{L:generalizedBS} are not enough in this situation.
One difference is that the laminations in \fullref{L:generalizedBS} were not canonical. 
We chose the laminations for convenience, there was no reason a quasi-isometry had to preserve lines of the lamination, but by making convenient choices we could build a quasi-isometry that did. In $B(R)$ we do not get to choose. 
The lamination is forced on us because we know the vertex spaces must be preserved. 

A second difference is that in \fullref{L:generalizedBS} the construction of quasi-isometries on the lines of the laminations was engineered so that we could get both height change and relative numbers of incident edges right. 
We can not hope to be so lucky in $B(R)$. For one thing there are adjacent vertices of different types that may have different edge densities. For another, each horizontal component represents a vertex space that belongs to a number of different P--sets. If we define the quasi-isometry to meet the needs of a particular P--set, we will not be able to glue up all of the different quasi-isometries at the end.
Because of these problems we will not be able to match adjacent lines in the laminations as we did in \fullref{L:generalizedBS}. 
However, it was not really necessary for a quasi-isometry to match adjacent lines of a lamination to adjacent lines of the other; it is good enough to take nearby lines to nearby lines.
This is the major difference between \fullref{L:qiofPsets} and \fullref{L:generalizedBS}, we get around the additional complications present in $B(R)$ essentially by showing that we can find a matching between suitable nearby horizontal components.

A third difference, which is a technical detail, is that in \fullref{L:generalizedBS}, long enough segments in the lines of the lamination always had strictly positive height change.
This was convenient because, when matching lines adjacent to a line $r$ to lines adjacent to a line $r'$, we did not have to worry about relative distances between lines adjacent to $r$ compared to the distances between their images adjacent to $r'$. Getting the heights right ensured that the distances would be right too. In $B(R)$ the lines of the lamination are called ``horizontal components'' exactly because there is no height change along them. This means we have to work a little harder to make sure the relevant distances work out correctly.

The next lemma gives a quasi-isometry of a P--set space for a P--set of bounded height change.
The various parameters just say that this is a quasi-isometry, built according to an extension $\mathcal{E}$, starting from a previously defined map on one of the depth zero vertex spaces.

\begin{lemma}\label{L:qiofPsets}
Let $G_1=\pi_1(\Gamma_1)$ and $G_2=\pi_1(\Gamma_2)$ be two tubular groups.
Let $R\subset D\varGamma_1$ and $S\subset D\varGamma_2$ be P--sets of bounded height change.
Let $x_0\in R$ and $y_0\in S$ be depth zero vertices.
Let $\mathcal{E}=(e_{ij})$ be an extension for $([R],[S])$.
Let $\alpha$ and $\beta$ be arbitrary positive real numbers.
Then there is a quasi-isometry $\Phi=\Phi^0_{\mathcal{E},(R,x_0,\alpha),(S,y_0,\beta)}\from X_R\to Y_S$ with the following properties:
\begin{enumerate}
\item $\Phi(X_{x_0})=Y_{y_0}$, and this map is, up to isometry, just a homothety with expansion by $\frac{\beta}{\alpha}$.
\item $\Phi$ induces a bijection, $\Phi_{\#}\from\mathcal{V}_0R\to \mathcal{V}_0S$.
\item For every depth zero vertex $v\neq x_0$, up to isometry, $\Phi |_{X_v}$ is a homothety with expansion by $\frac{\beta}{\alpha} 2^{h(v)-h(\Phi_{\#}(v))}$.
\item Excluding $x_0$ and $y_0$, there exist vertices of type $\{[R],i\}$ mapping to vertices of type $\{[S],j\}$ $\iff$ $e_{ij}$ is non-zero.
\item $\Phi$ is a bilipschitz bijection $(\mathcal{V}_0X_R,d_{X_R}|_{\mathcal{V}_0X_R})\to (\mathcal{V}_0Y_S,d_{Y_S}|_{\mathcal{V}_0Y_S})$.
\end{enumerate}
\end{lemma}
\begin{proof}
Pick some $J$ such that $\frac{J}{2}$ is a bound for the absolute value of the height change between any two vertices of $R$ or any two vertices of $S$.

For any $x\in R$, let $h(x)=h(x_0,x)$, and similarly for $y\in S$ with respect to $y_0$.

Let $x$ be a depth zero vertex of $R$.
Let $e$ be an edge of $R$ incident to $x$.
Since $\Stab_{G_1}(e)$ has infinite index in $\Stab_{G_1}(x)$, the orbit of $e$ by $\Stab_{G_1}(x)\cap \Stab_{G_1}(R)$ contains infinitely many other edges incident to $x$.
The edge spaces for these edges glue on to $X_x$ along parallel lines, and the distance between such lines can be bounded in terms of $[x]$, $[R]$, and $[e]$.

Pick some $L$ such that for each $P\in\{R, S\}$, and every $v\in \mathcal{V}_0P$, any open interval of length $L$ in $B(v)$ has at least two incident vertical edges from each $\Stab_G(v)\cap \Stab_G(P)$ equivalence class of edge incident to $v$, and at least three total incident vertical edges.

Choose coordinates for $X_R$ as discussed above. Choose coordinates for $Y_S$ in a similar fashion, with the following provision.
If we already have a map $X_{x_0}\to Y_{y_0}$, choose the coordinates of $Y_{y_0}$ so that the origin of $Y_{y_0}$ is the image of the origin of $X_{x_0}$ and so that map preserves the orientation of the second coordinate.

For every $v\in \mathcal{V}_0P$, do the following:
Define the basepoint $\xi_v$ of $B(v)$ to be the point of $B(v)$ with coordinate zero.
Slide any edge incident to the open interval \[\left(-\frac{1}{2}L\cdot 2^{h(v)},\frac{1}{2}L\cdot 2^{h(v)})\right)\] in $B(v)$ to $0$.
Note that for each $v$ this interval is of length $L$.

These sliding operations change distances between points in different horizontal components, but preserve distance within a fixed horizontal component.

For the vertices of the vertical components that are points of intersection with horizontal components, the \emph{type} of the vertex is just the type of the vertex in $D\varGamma$ corresponding to that horizontal component.
After sliding, the new vertical components are composed of unions of the original vertical components.
In particular, the new vertical components containing $\xi_{x_0}$ and $\xi_{y_0}$ are infinite, bounded valence trees; call them $V$ and $W$, respectively.

The image of $V$ in $\Gamma_1$ is the same as the image of $R$ in $\Gamma_1$.
Therefore, the set of vertices of $V$ of any particular vertex type are coarsely dense in $V$.
The diameter of the image of $R$ in $\Gamma_1$ provides a coarseness constant.

Similar statements are true for $W$ as well.
Let $\delta$ be the greater of the diameters of $R$ in $\Gamma_1$ and $S$ in $\Gamma_2$.

Every vertex of $V$ and $W$ is distance at most $\delta$ from a vertex of valence at least three, and there are no valence one vertices, so $V$ and $W$ are $\delta$-\emph{bushy}.
Also, the valences of $V$ and $W$ are bounded above.

For convenience we will assume that the extension $\mathcal{E}$ gives a bijection from vertex types of $R$ to vertex types of $S$.
To arrange this, note that the set of vertices of $V$ of a particular type can be subdivided into a finite number of subsets, each of which is also dense in $V$, and similarly for $W$.

Apply \fullref{L:3valent} to get a quasi-isometry $\phi \from V\to W$ which is a bilipschitz bijection on vertices and respects the partitions into vertex type.
The bilipschitz constant, $M_x$, depends only on the valence bounds, bushiness constants, density constants, and number of vertex types.

This process will be called extension along a vertical component.
Note that every depth zero vertex in $V$ or $W$ is the basepoint of its horizontal component.

Suppose extension along a vertical component identifies $\xi_x$ to $\xi_y$.

Let \[N_x=\max \{L\cdot 2^{h(x)},\frac{\alpha}{\beta}L\cdot 2^{h(y)}\}.\]
For $n=1,2,\dots $, consider the half-open intervals 
\begin{gather*}
\left[\frac{1}{2}L\cdot 2^{h(x)}+(n-1)N_x,\frac{1}{2}L\cdot 2^{h(x)}+nN_x \right) \subset B(x)\\
\left(-\frac{1}{2}L\cdot 2^{h(x)}-nN_x,-\frac{1}{2}L\cdot 2^{h(x)}-(n-1)N_x \right] \subset B(x)\\
\left[\frac{1}{2}L\cdot 2^{h(y)}+\frac{\beta}{\alpha}(n-1)N_x,\frac{1}{2}L\cdot 2^{h(y)}+\frac{\beta}{\alpha}nN_x \right) \subset B(y)\\
\left(-\frac{1}{2}L\cdot 2^{h(y)}-\frac{\beta}{\alpha}nN_x,-\frac{1}{2}L\cdot 2^{h(y)}-\frac{\beta}{\alpha}(n-1)N_x \right] \subset B(y)\\
\end{gather*}

Each of these intervals has length at least $L$, so there are at least three vertical edges incident to each interval.
On the other hand, the lengths of the intervals are bounded above by $N=L\cdot 2^J\max \{\frac{\alpha}{\beta},\frac{\beta}{\alpha}\}$.

In $B(x)$, slide all incident vertical edges to the closed endpoint of the interval to which they attach.
In $B(y)$, slide incident vertical edges in the $nth$ right open interval to the point $\frac{\beta}{\alpha}(\frac{1}{2}L\cdot 2^{h(x)}+(n-1)N_x)$.
Perform similar operations for the left open intervals of $B(y)$.

Extend $\phi$ by $\phi |_{B(x)}\from B(x)\to B(y)\from r \mapsto \pm \frac{\beta}{\alpha} r$.
The choice of orientation here is determined by the element $e_{ij}$ associated to $x,y$.
The edge sliding matches up the points to which vertical edges attach, and each of these will be the basepoint of a new vertical component to extend along.

Alternate extending along collections of vertical and horizontal components to build $\phi$.

We make the following observations about this construction:
\begin{enumerate}
\item Every vertical edge has at most one endpoint that slides.
\item No edge slides more than once.
\item Sliding an edge only changes distances between points separated by the edge.
\item No edge slides more than distance $N$.
\end{enumerate}

Therefore, distances in $B(\mathcal{V}_0R)$ and $B(\mathcal{V}_0S)$ change by at most a multiplicative factor of $1+N$.
Extension along horizontal components changes distance by at most a multiplicative factor of $\frac{\beta}{\alpha}2^J$.
Vertical components have valence bounded by $N$ and information from $X$ and $Y$, so there is a uniform bound $M$ on the bilipschitz constants $M_x$.
Thus, $\frac{\beta}{\alpha}(1+N)^2M2^J$ gives a bilipschitz constant for $\phi\from B(\mathcal{V}_0R)\to B(\mathcal{V}_0S)$.

Define $\Phi\from B(R)\times \mathbb{R} \to B(S)\times \mathbb{R}$ by $(t,u) \mapsto (\phi(t),\frac{\beta}{\alpha}u)$.

Geodesics in $X_R$ and $Y_S$ can be approximated within bounded multiplicative error by paths in which only one coordinate changes at a time, so $\Phi$ is bilipschitz on $\mathcal{V}_0X_R$.

The union of vertex spaces is dense in a P--set space, so $\Phi$ gives a quasi-isometry of P--set spaces.
\end{proof}

We have a similar result for P--sets of unbounded height change.
\begin{lemma}\label{L:qiofUnboundedPsets}
Let $G_1=\pi_1(\Gamma_1)$ and $G_2=\pi_1(\Gamma_2)$ be two tubular groups.
Let $R\subset D\varGamma_1$ and $S\subset D\varGamma_2$ be P--sets of unbounded height change.
Let $x_0\in R$ and $y_0\in S$ be depth zero vertices.
Let $\mathcal{E}=(e_{ij})$ be an extension for $([R],[S])$.
Let $\alpha$ and $\beta$ be arbitrary positive real numbers.
Then there is a quasi-isometry $\Phi=\Phi^{\infty}_{\mathcal{E},(R,x_0,\alpha),(S,y_0,\beta,)}\from X_R\to Y_S$ with the following properties:
\begin{enumerate}
\item $\Phi(X_{x_0})=Y_{y_0}$, and this map is, up to isometry, just a homothety with expansion by $\frac{\beta}{\alpha}$.
\item $\Phi$ induces a bijection, $\Phi_{\#}\from\mathcal{V}_0R\to \mathcal{V}_0S$.
\item Up to isometry, and for $v\neq x_0$, $\Phi |_{X_v}$ is a homothety with uniformly bounded expansion factor.
\item Excluding $x_0$ and $y_0$, there exist vertices of type $\{[R],i\}$ mapping to vertices of type $\{[S],j\}$ $\iff$ $e_{ij}$ is non-zero.
\item $\Phi$ is a bilipschitz bijection $(\mathcal{V}_0X_R,d_{X_R}|_{\mathcal{V}_0X_R})\to (\mathcal{V}_0Y_S,d_{Y_S}|_{\mathcal{V}_0Y_S})$.
\end{enumerate}
\end{lemma}
\begin{proof}
The proof differs from the proof of \fullref{L:qiofPsets} only in the how we extend along vertical components.

For $i=1,2$, consider the projection of the appropriate P--set to the graph $\Gamma_i$.
Let $\mu_i$ be the largest height change that occurs across a single edge in the projection.
The edge sliding is set up so that we can apply \fullref{L:generalizedBS} to find a coarsely height preserving quasi-isometry between the resulting vertical components.
Furthermore, the height error can be bounded in terms of the projections to the $\Gamma_i$.
Let $K$ be such a bound; $K$ plays the role in this case that the constant $J$ played in the bounded height change case.

Perform the edge sliding as in the previous Lemma.
Suppose $V$ and $W$ are vertical components based at $v_0\in B(x)$ and $w_0\in B(y)$, respectively.
Notice that the valence of $V$ at $v_0$ may depend on $\frac{\beta}{\alpha}$, since $v_0$ may not be equal to $\xi_x$.
However, every other vertex in $V\cap B(\mathcal{V}_0R)$ and $W\cap B(\mathcal{V}_0S)$ is the basepoint of its vertical component, so the valence of the trees at these points is bounded independently of $\frac{\beta}{\alpha}$.
It is possible to ``disperse'' the extra edges at $v_0$ by a height and type preserving bilipschitz bijection $V\cap B(\mathcal{V}_0R)$ to a new tree $V'$.
In the same way replace $W$ by a new tree $W'$.

There is some vertex $\xi_w$ in $W$ at bounded distance from $\xi_{w_0}$ such that 
\[|h(\xi_{v_0},\xi_{w_0})+\log_2\frac{\beta}{\alpha}|\leq \mu_2\]
Let $\phi$ be a $(K+\mu_2)$-coarsely height preserving quasi-isometry between $V'$ and $W'$ that maps $v_0$ to $\xi_w$.
As in the previous Lemma, $\phi$, can then be changed a bounded amount to give a bijection $\phi'$ of depth zero vertices that respects the partitions by type.
Moreover, the amount $\phi$ needs to be changed is independent of $\frac{\beta}{\alpha}$, so the height error is still independent of $\frac{\beta}{\alpha}$.

The map $V\to V' \to W' \to W$ is then a bilipschitz bijection:
\[V\cap B(\mathcal{V}_0R) \to W\cap B(\mathcal{V}_0S)\]
Finally, we adjust this map to ensure $v_0$ maps to $w_0$.
The result is the map we will use to extend along vertical components.
Note that while the bilipschitz constants depend on $\frac{\beta}{\alpha}$, the height error at points other than $v_0$ is bounded independently of $\frac{\beta}{\alpha}$.
\end{proof}

\begin{remark}
Note the important difference in these two Lemmas. 
In both cases we start by identifying $x_0$ to $y_0$.
The parameters $\alpha$ and $\beta$ provide an initial height error; think of $E=\log_2(\frac{\beta}{\alpha})$ as the height error at $x_0$.

In the bounded case, the height error of the other vertices is in $[E-J, E+J]$.

In the unbounded case, $E$ influences the quasi-isometry constants, but the height error of the other vertices is independent of $E$. 
\end{remark}

\subsubsection{Putting the Pieces Together}\label{ss:equiv}
The following Theorem shows that a consistent set of strategies can be used to build a quasi-isometry between tubular groups.
\fullref{L:qiofPsets} and \fullref{L:qiofUnboundedPsets} provide the basic building blocks, and the quasi-isometry constants are controlled by controlling height error.
Two-line, depth zero vertices and P--sets of unbounded height change are flexible enough that height error will not be an issue.
In these cases we can choose the quasi-isometry to immediately ``correct'' any accumulated height error back to some uniformly bounded amount.
The danger of compounding height error comes from P--sets of bounded height change connected by depth zero vertices with at least three lines.
These height errors are controlled by the set of strategies.

\begin{theorem}\label{T:Construction}
For $i=1,2$, let $G_i=\pi_1(\Gamma_i)$ be a tubular group and $T_i$ its tree of P--sets.
If there is a consistent set of strategies for $T_1$ and $T_2$, then $G_1$ and $G_2$ are quasi-isometric.
\end{theorem}
\begin{proof}
Let $\mathcal{G}$ be a consistent set of strategies for $T_1$ and $T_2$.
We build a quasi-isometry $\Phi\from X\to Y$ and an allowable tree isomorphism $\phi\from T_1\to T_2$.

Let $\mathcal{M}_i$ be the matches of P--sets of bounded height change occurring in $\mathcal{G}$, and let $\{U_i\}$ and $\{L_i\}$ be height error bounds from the consistency check.

Choose a uniform $K'$ such that, for any P--sets $R\in T_1$ and $S\in T_2$ of unbounded height change, there is a quasi-isometry $\Phi^\infty\from X_R\to Y_S$ as in \fullref{L:qiofUnboundedPsets} with height error bounded by $K'$.

In $D\varGamma_1$ and $D\varGamma_2$ there are finitely many equivalence classes of P--sets of bounded height change, so there is some maximum height change, $K''$, that can occur between depth zero vertices in such a P--set. This implies that for any P--set $R\in T_1$ of bounded height change, any depth zero vertices $v$, $v'\in R$ and any isomorphism $\phi\from T_1\to T_2$ we have $|\Err_\phi(v)-\Err_\phi(v')|\leq 2K''$ and $|\Err_\phi(v')-\Err_\phi(R)|\leq K''$, where $\Err_\phi(v)=h(v)-h(\phi(v))$.

Let $K=\max\{K',K''\}$. Assume the $U_i$ are greater than $2K$ and the $L_i$ are less than $-2K$.

We build a quasi-isometry inductively P--set space by P--set space.

There are a number of cases to consider, but for all cases the outline of the construction is the same. In the cases that follow, we will refer to this process as ``extending the maps''.

At each step we are given P--sets $R,\,R'\in T_1$, and a depth zero vertex $v\in R\cap R'$. 
We assume that we have a partial tree isomorphism $\phi\from T_1\to T_2$ that is defined at least on $R$, $R'$ and $v$, and a quasi-isometry $\Phi_R\from X_R\to Y_{\phi(R)}$.
We want to extend $\phi$ to a neighborhood of $R'$ and define a quasi-isometry $\Phi_{R'}\from X_{R'}\to Y_{\phi(R')}$.

Define a quasi-isometry $\Phi_{R'}\from X_{R'}\to Y_{\phi(R')}$ by consulting $\mathcal{G}$ for an ``appropriate'' strategy, $\mathcal{S}$, for the match $([R'],[\phi(R')])$.
If $R'$ is of unbounded height change there is only one strategy to choose. If $R'$ is of bounded height change then we choose the positive strategy if $\Err_\phi(R')<0$ and the negative strategy if $\Err_\phi(R')\geq0$.
Associated to $\mathcal{S}$ is an extension $\mathcal{E}$. 
Define $\Phi_{R'}=\Phi^0_{\mathcal{E},(R',v,\alpha),(\phi(R'),\phi(v),\beta)}$ from \fullref{L:qiofPsets} if $R'$ is of bounded height change or $\Phi_{R'}=\Phi^\infty_{\mathcal{E},(R',v,\alpha),(\phi(R'),\phi(v),\beta)}$ from \fullref{L:qiofUnboundedPsets} if $R'$ is of unbounded height change.
The constants $\alpha$ and $\beta$ are to be determined in the induction steps.

Extend $\phi$ from $R'$ according to $\mathcal{E}$, using the bijection of depth zero vertices furnished by $(\Phi_{R'})_{\#}$.

 In each case below we check that if we started with a P--set $R'$ such that $([R'],[\phi(R')])=\mathcal{M}_i$ then all depth zero vertices of $R'$ have height error in the interval $[L_i-K, U_i+K]$ if $R'$ is of bounded height change or $[-K,K]$ if $R'$ is of unbounded height change.

\case{Base Case }{
Pick any P--set vertices $R\in T_{1}$ and $S\in T_2$ such that $([R],[S])$ is the match of a vertex of $\mathcal{G}$.
Set $\phi(R)=S$.

Pick any depth zero vertex $v\in R$.
Declare $v$ to be the basepoint, so $v$ has height zero and all other heights are determined relative to $v$.

Consult $\mathcal{G}$ for an appropriate strategy $\mathcal{S}$ for $([R],[S])$ with extension $\mathcal{E}$.
If $v$ is of type $\{[R],i\}$ with respect to $R$, find some $j$ such that the $i,j$ entry of $\mathcal{E}$ is nonzero and choose a depth zero vertex $w\in S$ of type $\{[S],j\}$.
Define $\phi(v)=w$.

Declare $w$ to be height zero, and define heights in $S$ relative to $w$.

Now we have a P--set $R$, a preferred vertex $v\in R$ and a partial tree isomorphism defined at least for $v$ and $R$, with $\Err_\phi(v)=0$.
Extend the maps with $\alpha=\beta=1$.

If $R$ is of unbounded height change then by construction the height errors for all the depth zero vertices of $R$ will be in the interval $[-K,K]$.
If $R$ is of bounded height change then for any depth zero vertex $v'\in R$, 
\[|\Err_\phi(v)-\Err_\phi(v')|\leq 2K\implies \Err_\phi(v')\in[-2K,2K]\] 
and \[|\Err_\phi(v)-\Err_\phi(R)|\leq K\implies \Err_\phi(R)\in[-K,K]\]
}

Now there are a number of situations to consider, depending on whether $R$ and $R'$ are of bounded or unbounded height change, and whether $v$ has two lines or more than two lines.

\case{$v$ has 2 lines}{
Define height in $R'$ relative to $v$ and define height in $\phi(R')$ relative to $\phi(v)$.

Extend the maps with $\alpha=\beta=1$.

The height errors of depth zero vertices of $R'$ are in $[-2K,2K]$.

Furthermore, if $R'$ is of bounded height change then $\Err_\phi(R')\in [-K,K]$.
}

In the remaining cases, $v$ has more than two lines, so $R$ and $R'$ are in the same rigid component.
Height in $R'$ is already defined.
The same goes for $\phi(R)$ and $\phi(R')$.

\case{$R'$ of unbounded height change}{
Extend the maps with $\alpha=2^{-h(v)}$ and $\beta=2^{-h(\phi(v))}$.

By construction, the height errors of depth zero vertices of $R'$ (other than $v$) are in $[-K,K]$.
}

\case{$R$ of unbounded height change, $R'$ of bounded height change}{
Since $R$ was of unbounded height change, we can assume that the height error of $v$ is in $[-K,K]$.

Extend the maps with $\alpha=2^{-h(v)}$ and $\beta=2^{-h(\phi(v))}$.

The height error of $R'$ is in $[-2K,2K]$.

The height errors of depth zero vertices of $R'$ are in $[-3K,3K]$.
}

In every case thus far we have produced a universal bound for the height errors of depth zero vertices. 
The final case is the one to worry about; this is the case where height error can compound.

\case{$R$ and $R'$ of bounded height change}{
Suppose $([R],[\phi(R)])=\mathcal{M}_i$, and suppose $([R'],[\phi(R')])=\mathcal{M}_j$.
Assume the height error of $R$ was in $[L_i,U_i]$.

Extend the maps with $\alpha=2^{-h(v)}$ and $\beta=2^{-h(\phi(v))}$.

A priori, if $v'$ is a depth zero vertex of $R'$ we know
$|\Err_\phi(v')-\Err_\phi(R')|\leq K$ and $|\Err_\phi(R)-\Err_\phi(R')|\leq 2K$, so the height errors could be worse than in the previous step.
If the errors continue to get worse in each step then we lose control.

However, because we have built according to the set of strategies, we are guaranteed that the height error of $R'$ is in $[L_j,U_j]$, 
so $\Err_\phi(v')\in [L_j-K,U_j+K]$.
}

This completes the induction steps.

The various maps $\Phi_R$ were constructed so that they agree when they overlap on any depth zero vertex spaces of at least three lines.
They may not agree on two-line vertex spaces, but this can be fixed, since two-line patterns are not rigid.

Suppose $v$ is a two-line vertex joining P--sets $R$ and $R'$.
Suppose $\Phi_R|_{X_v}$ is expansion by a factor $C$, and $\Phi_{R'}|_{X_v}$ is expansion by a factor $C'$.
Change these maps so that they both expand by a factor of $C$ along the lines where the $R'$ edges glue on, and expand by a factor of $C'$ along the lines where the $R$ edges glue on.
This does not change the maps that we had on $B(v)\subset B(R)$ and $B(v)\subset B(R')$.

The map $\Phi$ is built from various pieces $\Phi_R$.
The $\Phi_R$ are bilipschitz bijections on $\mathcal{V}_0X_R$, and agree on vertex spaces where their domains intersect.
We have arranged that the bilipschitz constants of $\Phi_R$ depend on $[R]$ and $[\phi(R)]$ and on height error.
Tubular groups come from finite graphs of groups, so there are only finitely many equivalence classes of P--set.
Thus, invariants of equivalence classes of P--sets can be uniformly bounded.
Height error is bounded by $K+\max_i \{U_i, |L_i|\}$.

The map $\Phi$ is thus a bilipschitz bijection $\mathcal{V}_0X\to \mathcal{V}_0Y$.
The set $\mathcal{V}_0X$ is dense in $X$, and $\mathcal{V}_0Y$ is dense in $Y$, so $\Phi\from X\to Y$ is a quasi-isometry.
\end{proof}

\begin{maintheorem}
For $i=1,2$, let $G_i=\pi_1(\Gamma_i)$ be tubular groups.
There is an algorithm that in finite time determines whether or not $G_1$ and $G_2$ are quasi-isometric.
\end{maintheorem}
\begin{proof}
Run the algorithm from \fullref{T:Algorithm}.
This algorithm halts in finite time.
If the groups are quasi-isometric then \fullref{Corollary:qiimpliesallowable} and \fullref{L:Success} guarantee that the algorithm succeeds and produces a consistent set of strategies.
Conversely, \fullref{T:Construction} builds a quasi-isometry from a consistent set of strategies, so the algorithm succeeds only if the groups are quasi-isometric.
\end{proof}

\section{Consequences and Examples}\label{C:CandE}

\begin{example}\label{Ex:RAAG}
Consider two tubular groups all of whose depth zero vertices have two lines.
The only obstruction to building a quasi-isometry in this case is that P--sets of bounded height change must be matched to P--sets of bounded height change, and similarly for P--sets of unbounded height change.
\end{example}

In particular, we recover the following result of Behrstock and Neumann \cite{BehNeu06}:
\begin{corollary}\label{Cor:RAAG}
 Any two Right Angled Artin Groups whose defining graphs are trees of diameter at least three are quasi-isometric.
\end{corollary}
This follows since such groups have graph of groups decompositions as tubular groups where all depth zero vertices have two lines and all P--sets are of bounded height change.

Suppose $G$ is a tubular group with at least three lines in every vertex and having bounded height change in every P--set.
Then there is a notion of height change for every edge in the tree of P--sets.

If $\gamma$ is a geodesic ray in $T_G$, $\gamma$ \emph{has coarse slope} $m$, $\slope (\gamma)=m$, if there exists a $C>0$ such that for all $t$, \[tm-C\leq h(\gamma(0),\gamma(t)) \leq tm+C\]
\begin{proposition}
There is some $m_G=\max_{\gamma\subset T_G}\slope (\gamma)$, and the set of slopes of geodesic rays of $T_G$ is dense in the interval $[-m_G, m_G]$.
\end{proposition}
\begin{proof}
Let $\gamma$ be a geodesic ray in $T_G$.
Consider the image of $\gamma$ in the graph that is the quotient of $T_G$ by the group action.
This graph is finite, so the image of $\gamma$ can be written as some initial segment of bounded length followed by loops in the graph.
The initial segment contributes a bounded amount to height, so it can be discarded.

There are only finitely many simple closed curves in a finite graph, so pick one, $\alpha$, of maximal slope.
It is not hard to see that any loop can be reduced to a simple closed curve without reducing slope.
Therefore, $\slope (\gamma)\leq \slope (\alpha)$.
Then $m_G=\slope (\alpha)$.

Loops that are homotopically trivial have zero height change, so we can get a geodesic ray in $T_G$ of slope arbitrarily close to any $m\in [-m_G,\, m_G]$ by concatenating lifts of $\pm \alpha$ with lifts of trivial loops. 
\end{proof}

\begin{corollary}
There is a maximum coarse slope in each rigid component of $T_G$ that contains no P--set of unbounded height change.
\end{corollary}

Suppose $G_1$ and $G_2$ are tubular groups with at least three lines in every vertex and having bounded height change in every P--set.
Also, suppose that $G_1$ and $G_2$ have the same sets of equivalence classes of edge patterns, and that each of these patterns has a unique (up to isometry) symmetric representative.
This occurs, for instance, if every vertex has three or four lines.

Suppose $m_{G_1}>m_{G_2}$, and $\gamma$ is a geodesic ray of slope $m_{G_1}$ in $T_{G_1}$.
Under a coarsely height preserving tree isomorphism, the image of $\gamma$ will still be a geodesic ray of slope $m_{G_1}$, contradicting maximality of $m_{G_2}$.
This means $m_G$ is a quasi-isometry invariant of $G$.

When vertices of $G$ have edge patterns without a unique symmetric representative, the value of $m_G$ depends on the choice of representatives.

The max slope is not a complete invariant, even among groups from the family described in Section~\ref{SSS:existence}.
In \fullref{Ex:bothslopes} we classify this family completely, and see that groups with the same max slope need not be quasi-isometric.

\begin{note}
In diagrams of trees of P--sets, P--sets are represented by large open circles labeled by a representative of the equivalence class of P--set.
One adjacent vertex is added for each P--set-stabilizer equivalence class of depth zero vertex contained in the P--set, and these are labeled by vertex types with respect to the more heavily shaded P--set in the center of the diagram. 
\end{note}

\begin{example}\label{Ex:bothslopes}
Suppose $G_1$ and $G_2$ are groups as in \fullref{Fi:OneTorusStandard}.
Each $G_i$ has only one equivalence class of P--set, $[C_i]$, and each P--set has three vertex types.

\begin{figure}[ht!]
\labellist
\small
\pinlabel $\{[C_i],2\}$ [tl] at 241 585
\pinlabel $\{[C_i],3\}$ [tr] at 316 570
\pinlabel $\{[C_i],1\}$ [lb] at 341 615
\pinlabel $C_i$ at 316 616
\pinlabel $C_i$ at 216 584
\pinlabel $C_i$ at 240.1 628.9
\pinlabel $C_i$ at 340.1 571
\pinlabel $C_i$ at 340.1 603
\pinlabel $C_i$ at 316 647.5
\pinlabel $C_i$ at 340 660.1
\pinlabel $\lambda_i$ [rB] at 200 659
\pinlabel $\lambda_i-\mu_i$ [rB] at 200 647
\pinlabel $\mu_i$ [rB] at 200 629
\pinlabel $0\phantom{_i}$ [rB] at 200 614
\pinlabel $-\mu_i$ [rB] at 200 602
\pinlabel $-\lambda_i+\mu_i$ [rB] at 200 583
\pinlabel $-\lambda_i$ [rB] at 200 570
\endlabellist
\centering
  \includegraphics{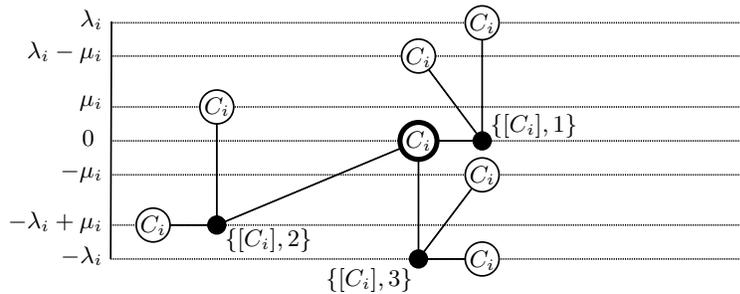}
\caption{Vertex types and relative heights for $[C_i]$.}
\label{Fi:classCitypes}
\end{figure}

First assume that $\lambda_i\geq \mu_i \geq 0$.

\fullref{Fi:classCitypes} summarizes the setup. Notice that the height change from a vertex of type $\{[C_i],3\}$ to a vertex of type $\{[C_i],2\}$ is $\mu_i$, and the height change from a vertex of type $\{[C_i],3\}$ to a vertex of type $\{[C_i],1\}$ is $\lambda_i$.

The maximum slope in $T_i$ is $\lambda_i$, so assume that $\lambda=\lambda_1=\lambda_2$.

If $\lambda=0$ then there are no height changes to worry about, so the groups are quasi-isometric, so assume $\lambda>0$.

First we attempt to build an extension $\mathcal{E}=(e_{jk})$ that gives a positive strategy.

To build a consistent set of strategies we will need a positive strategy that gives all non-negative height errors.
The largest negative height change that occurs is $-\lambda$ in both trees, so we must identify these.
To accomplish this, $e_{33}$ must be non-zero.
Furthermore, if we choose this entry to identify the P--sets at height $-\lambda$, we are forced to identify P--sets at height $-\mu_1$ with P--sets of height $-\mu_2$.
Thus, we need $\mu_1\leq \mu_2$.

We get the reverse inequality by considering a negative strategy.
We need a negative strategy that gives all non-positive height errors.
The largest positive height change that occurs is $\lambda$ in both trees, so we must identify these.
This forces an identification of P--sets at height $\lambda-\mu_1$ with P--sets of height $\lambda-\mu_2$, so we need $\mu_1\geq \mu_2$.

Thus, if $\mu_1\neq \mu_2$, the two groups can not be quasi-isometric.

When $\lambda_1=\lambda_2$ and $\mu_1=\mu_2$ the obvious strategy produces zero height error, so the groups are quasi-isometric.

We started by assuming that $\lambda$ and $\mu$ were non-negative.
It is not hard to see that groups from this family having pairs of height changes $(\lambda, \mu)$, $(\lambda-\mu, -\mu)$ and $(-\lambda, \mu-\lambda)$ are all quasi-isometric.
There is only one distinct such pair with both entries non-negative.
Therefore, the (unordered) pair of non-negative height changes gives a complete quasi-isometry invariant for this family.
Notice that the max slope is the larger of the two height changes, so we have many examples of groups with the same max slope that are not quasi-isometric.
\end{example}

Finally, we give an example of quasi-isometric groups where there are actually bounded height errors (not always 0 as in \fullref{Ex:bothslopes}) and we really need two distinct strategies for some match.
\begin{example}\label{Ex:noncommensurable}
Let $G_1=\pi_1(\Gamma_1)$ be a one torus group as in \fullref{Ex:bothslopes} with height changes $\lambda>\mu>0$ across the two edges.

Let $G_2=\pi_1(\Gamma_2)$ be a tubular group whose underlying graph is shown in \fullref{Fi:TwoTorus}. In this diagram the labels are height changes and an arc joining edges incident to a common vertex indicates that the corresponding edge strips glue to the vertex space along parallel lines.

\begin{figure}[ht!]
\labellist
\small \hair 2pt
\pinlabel $\mu$ [rt] at 99 717
\pinlabel $\lambda -\mu$ [lt] at 284 717
\pinlabel $\lambda$ [rt] at 189 687
\pinlabel $\lambda$ [bl] at 190 752
\endlabellist
\centering
  \includegraphics{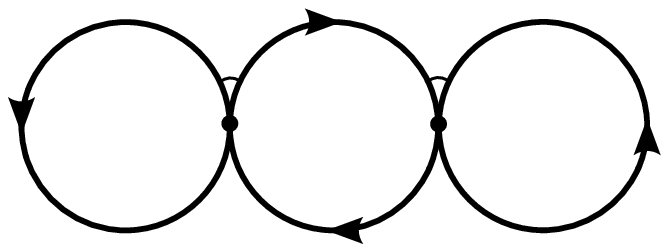}
\caption{A graph for a two torus group.}
\label{Fi:TwoTorus}
\end{figure}

The Bass-Serre tree, $D\varGamma_2$, has two equivalence classes of P-sets.
One class, $[A]$, consists of P-sets that project in $\Gamma_2$ to the union of the three edges connected by arcs.
The other, $[B]$, consists of P-sets that project to the remaining edge in $\Gamma_2$.

\fullref{Fi:classAtypes} and \fullref{Fi:classBtypes} define the vertex types for vertices in $A$ and $B$, respectively.

\begin{figure}[ht!]
\labellist
\small
\pinlabel $\{[A],1\}$ [lb] at 341 615
\pinlabel $\{[A],2\}$ [tl] at 241 585
\pinlabel $\{[A],3\}$ [tr] at 391 585
\pinlabel $\{[A],4\}$ [tr] at 316 570
\pinlabel $A$ at 316.5 616.3
\pinlabel $A$ at 316.5 648.5
\pinlabel $A$ at 340.5 603.8
\pinlabel $A$ at 216.5 584.8
\pinlabel $A$ at 392.5 629.3
\pinlabel $B$ at 340.3 571.5
\pinlabel $B$ at 240.3 629
\pinlabel $B$ at 415.8 584.5
\pinlabel $B$ at 340.3 660.5
\pinlabel $\lambda$ [rB] at 200 659
\pinlabel $\lambda-\mu$ [rB] at 200 647
\pinlabel $\mu$ [rB] at 200 628
\pinlabel $0$ [rB] at 200 614
\pinlabel $-\mu$ [rB] at 200 602
\pinlabel $-\lambda+\mu$ [rB] at 200 583
\pinlabel $-\lambda$ [rB] at 200 570
\endlabellist
\centering
  \includegraphics{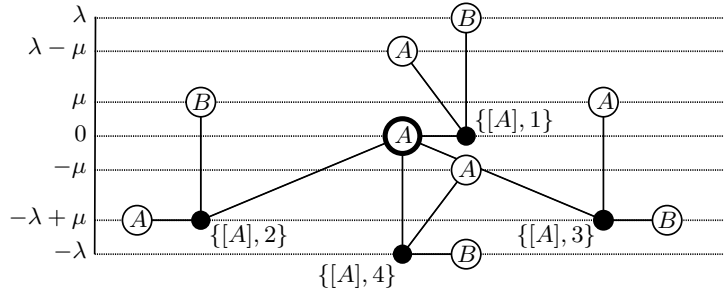}
\caption{Vertex types and relative heights for $[A]$.}
\label{Fi:classAtypes}
\end{figure}

\begin{figure}[ht!]
\labellist
\small
\pinlabel $\{[B],1\}$ [lb] at 341 615
\pinlabel $\{[B],2\}$ [tr] at 316 570
\pinlabel $B$ at 316 616
\pinlabel $A$ at 316.5 648.5
\pinlabel $A$ at 340.5 603.5
\pinlabel $A$ at 340.5 571.8
\pinlabel $A$ at 340.5 660.8
\pinlabel $\lambda$ [rB] at 200 659
\pinlabel $\lambda-\mu$ [rB] at 200 647
\pinlabel $\mu$ [rB] at 200 628
\pinlabel $0$ [rB] at 200 614
\pinlabel $-\mu$ [rB] at 200 602
\pinlabel $-\lambda+\mu$ [rB] at 200 583
\pinlabel $-\lambda$ [rB] at 200 570
\endlabellist
\centering
  \includegraphics{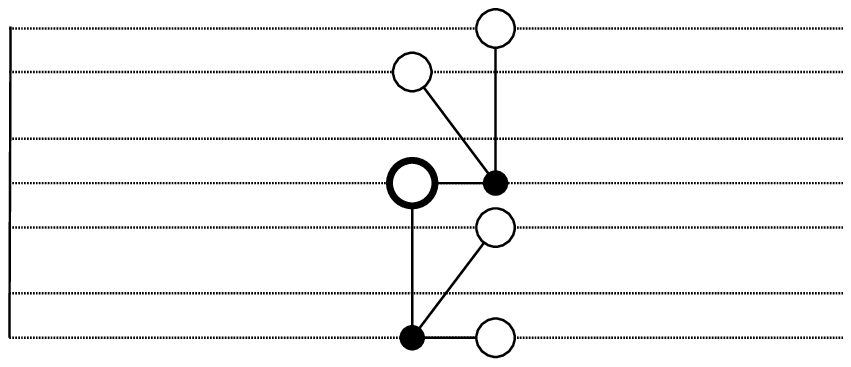}
\caption{Vertex types and relative heights for $[B]$.}
\label{Fi:classBtypes}
\end{figure}

For $([C],[A])$, consider the extension $\mathcal{E}_A=(e_{ij})$ defined as follows:
\[
\begin{matrix}
e_{11}=\left\{\begin{matrix}
\{[C],1\}\mapsto \{[A],1\}\\
\{[C],2\}\mapsto \{[A],2\}\\
\{[C],3\}\mapsto \{[B],2\}
\end{matrix}\right\}
&
e_{22}=\left\{\begin{matrix}
\{[C],1\}\mapsto \{[A],1\}\\
\{[C],2\}\mapsto \{[A],2\}\\
\{[C],3\}\mapsto \{[B],2\}
\end{matrix}\right\}
\\
e_{23}=\left\{\begin{matrix}
\{[C],1\}\mapsto \{[B],2\}\\
\{[C],2\}\mapsto \{[A],3\}\\
\{[C],3\}\mapsto \{[A],4\}
\end{matrix}\right\}
&
e_{34}=\left\{\begin{matrix}
\{[C],1\}\mapsto \{[B],1\}\\
\{[C],2\}\mapsto \{[A],3\}\\
\{[C],3\}\mapsto \{[A],4\}
\end{matrix}\right\}
\\
e_{ij}=\phantom{\lgroup}\begin{matrix}
\phantom{\{[C],1\}\mapsto \{[B],1\}}\\
0 \text{ for all other }ij\\
\\
\end{matrix}
\phantom{\lgroup}
\end{matrix}
\]

This extension says that given P--sets $C$ and $A$, we should split the vertices of type $\{[C],2\}$ with respect to $C$ and map half of them to the vertices of type $\{[A],2\}$ with respect to $A$ and the other half to vertices of type $\{[A],3\}$ with respect to $A$.

\fullref{Fi:inducedmappingA} shows the result of extending by $\mathcal{E}_A$.

\begin{figure}[ht!]
\labellist
\small
\pinlabel $C$ at 315.5 616.3
\pinlabel $C$ at 315.5 648.2
\pinlabel $C$ at 339.5 603.5
\pinlabel $C$ at 215.5 584.5
\pinlabel $C$ at 391.2 629
\pinlabel $C$ at 339.5 571.5
\pinlabel $C$ at 239.5 629
\pinlabel $C$ at 415.3 584.5
\pinlabel $C$ at 339.5 660.5
\pinlabel $A$ at 328.8 616.3
\pinlabel $A$ at 328.8 648.2
\pinlabel $A$ at 352.8 603.5
\pinlabel $A$ at 228.8 584.5
\pinlabel $A$ at 404.8 629
\pinlabel $B$ at 353 571.5
\pinlabel $B$ at 253 629
\pinlabel $B$ at 428.5 584.5
\pinlabel $B$ at 353 660.5
\pinlabel $\lambda$ [rB] at 200 659
\pinlabel $\lambda-\mu$ [rB] at 200 647
\pinlabel $\mu$ [rB] at 200 628
\pinlabel $0$ [rB] at 200 614
\pinlabel $-\mu$ [rB] at 200 602
\pinlabel $-\lambda+\mu$ [rB] at 200 583
\pinlabel $-\lambda$ [rB] at 200 570
\endlabellist
\centering
\includegraphics{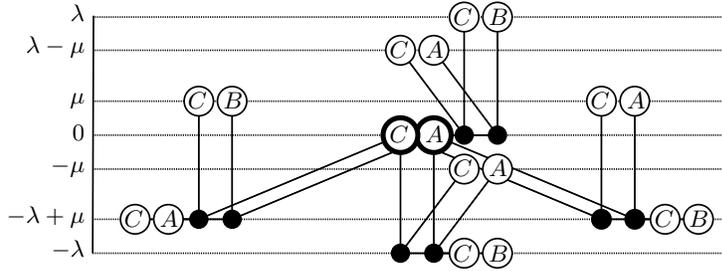}
\caption{Extension by $\mathcal{E}_A$.}
\label{Fi:inducedmappingA}
\end{figure}

All the matches induced by this extension have height error 0, so the strategy $\mathcal{S}_A$ for $\mathcal{E}_A$ can be used as a positive and a negative strategy.

For $([C],[B])$ we define extensions $\mathcal{E}_B^+$ and $\mathcal{E}_B^-$, extension by which are depicted in \fullref{Fi:inducedmappingBplus} and \fullref{Fi:inducedmappingBminus}, respectively.

\begin{figure}[ht!]
\labellist
\small
\pinlabel $C$ at 315.5 616.3
\pinlabel $C$ at 315.5 648.2
\pinlabel $C$ at 339.5 603.5
\pinlabel $C$ at 215.5 584.5
\pinlabel $C$ at 339.5 571.5
\pinlabel $C$ at 239.5 629
\pinlabel $C$ at 339.5 660.5
\pinlabel $B$ at 328.8 616.3
\pinlabel $A$ at 328.8 648.2
\pinlabel $A$ at 352.8 603.5
\pinlabel $A$ at 240.5 571.5
\pinlabel $A$ at 353 571.5
\pinlabel $A$ at 264 603.5
\pinlabel $A$ at 353 660.5
\pinlabel $\lambda$ [rB] at 200 659
\pinlabel $\lambda-\mu$ [rB] at 200 647
\pinlabel $\mu$ [rB] at 200 628
\pinlabel $0$ [rB] at 200 614
\pinlabel $-\mu$ [rB] at 200 602
\pinlabel $-\lambda+\mu$ [rB] at 200 583
\pinlabel $-\lambda$ [rB] at 200 570
\endlabellist
\centering
\includegraphics{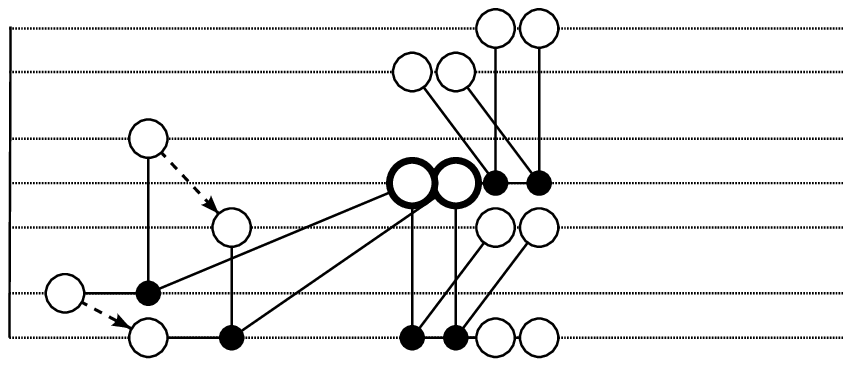}
\caption{Extension by $\mathcal{E}_B^+$.}
\label{Fi:inducedmappingBplus}
\end{figure}

The extension $\mathcal{E}_B^+$ induces matches $([C],[A])$ with height errors 0, $\mu$, and $2\mu$, so it gives rise to a positive strategy.

\begin{figure}[ht!]
\labellist
\small
\pinlabel $C$ at 315.5 616.3
\pinlabel $C$ at 315.5 648.2
\pinlabel $C$ at 339.5 603.5
\pinlabel $C$ at 215.5 584.5
\pinlabel $C$ at 339.5 571.5
\pinlabel $C$ at 239.5 629
\pinlabel $C$ at 339.5 660.5
\pinlabel $B$ at 328.8 616.3
\pinlabel $A$ at 328.8 648.2
\pinlabel $A$ at 352.8 603.5
\pinlabel $A$ at 240.5 648.2
\pinlabel $A$ at 353 571.5
\pinlabel $A$ at 264 660.5
\pinlabel $A$ at 353 660.5
\pinlabel $\lambda$ [rB] at 200 659
\pinlabel $\lambda-\mu$ [rB] at 200 647
\pinlabel $\mu$ [rB] at 200 628
\pinlabel $0$ [rB] at 200 614
\pinlabel $-\mu$ [rB] at 200 602
\pinlabel $-\lambda+\mu$ [rB] at 200 583
\pinlabel $-\lambda$ [rB] at 200 570
\endlabellist
\centering
\includegraphics{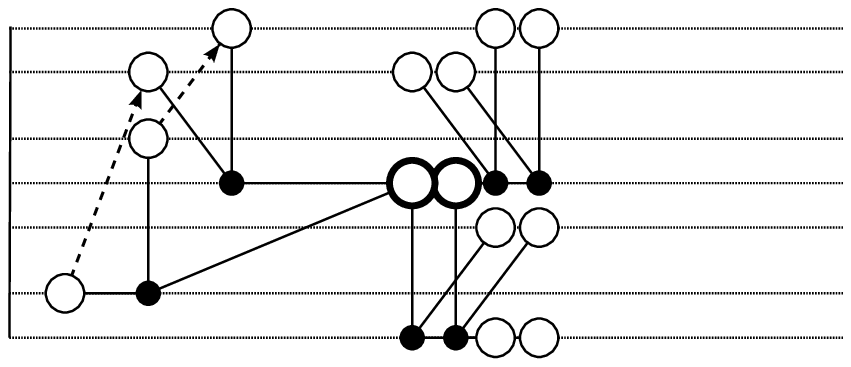}
\caption{Extension by $\mathcal{E}_B^-$.}
\label{Fi:inducedmappingBminus}
\end{figure}

The extension $\mathcal{E}_B^-$ induces matches $([C],[A])$ with height errors 0, $-2(\lambda-\mu)$, $\lambda-\mu$, so it gives rise to a negative strategy.

These strategies give the set of strategies in \fullref{Fi:setofstrat}.

\begin{figure}[ht!]
\labellist
\small
\pinlabel $([C],[A])$ [r] at 286 325
\pinlabel $([C],[B])$ [l] at 425 325
\pinlabel $\pm$ [r] at 207 325
\pinlabel $0$ [l] at 208 325
\pinlabel $+$ [b] at 356 378
\pinlabel $2\mu$ [t] at 356 378
\pinlabel $\pm$ [b] at 356 325
\pinlabel $0$ [t] at 356 325
\pinlabel $-$ [b] at 356 274
\pinlabel $-2(\lambda-\mu)$ [t] at 356 274
\endlabellist
\centering
\includegraphics{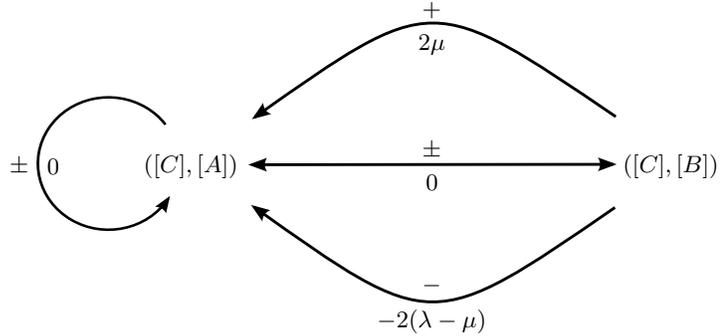}
\caption{The set of strategies}
\label{Fi:setofstrat}
\end{figure}

This set of strategies gives rise to a system of inequalities:
\[
\begin{gathered}
L_1 \leq M_1 \leq U_1\\
L_2 \leq M_2 \leq U_2\\
U_1 + 0 \leq U_2\\
M_1 + 0 \geq L_2\\
\end{gathered}
\qquad
\begin{gathered}
M_1 + 0 \leq U_2\\
L_1 + 0 \geq L_2\\
M_2+2\mu \leq U_1\\
L_2+2\mu \geq L_1
\end{gathered}
\qquad
\begin{gathered}
U_2 - 0 \leq U_1\\
M_2 - 0 \geq L_1\\
U_2-2(\lambda-\mu) \leq U_1\\
M_2-2(\lambda-\mu) \geq L_1
\end{gathered}
\]

This system has solutions, for instance:
\[\begin{gathered}
U_1=U_2=2\lambda\\
M_1=M_2=0\\
L_1=L_2=-2\lambda
\end{gathered}\]
Therefore, groups of this form are quasi-isometric.
\end{example}

\bibliographystyle{hyperamsplain}
\bibliography{masterbib}

\providecommand{\bysame}{\leavevmode\hbox to3em{\hrulefill}\thinspace}
\providecommand{\MR}{\relax\ifhmode\unskip\space\fi MR }
\providecommand{\MRhref}[2]{%
  \href{http://www.ams.org/mathscinet-getitem?mr=#1}{#2}
}
\providecommand{\href}[2]{#2}
\begin{thebibliography}{10}

\bibitem{BehNeu06}
Jason~A. Behrstock and Walter~D. Neumann, \emph{Quasi-isometric classification
  of graph manifold groups}, Duke Math. J. \textbf{141} (2008), 217--240,
  \href{http://arXiv.org/abs/math.GT/0604042}{\texttt{arXiv:math.GT/0604042}}.

\bibitem{BraBri00}
N.~Brady and M.~R. Bridson, \emph{There is only one gap in the isoperimetric
  spectrum}, Geom. Funct. Anal. \textbf{10} (2000), no.~5, 1053--1070.

\bibitem{BraBriFor09}
Noel Brady, Martin Bridson, Max Forester, and Krishnan Shankar, \emph{Snowflake
  groups, {Perron-Frobenius} eigenvalues, and isoperimetric spectra}, Geometry
  \& Topology \textbf{13} (2009), to appear,
  \href{http://arXiv.org/abs/math.GR/0608155}{\texttt{arXiv:math.GR/0608155}}.

\bibitem{BriGer96}
M.~R. Bridson and S.~M. Gersten, \emph{The optimal isoperimetric inequality for
  torus bundles over the circle}, Q. J. Math. \textbf{2} (1996), no.~47, 1--23.

\bibitem{BriHae99}
Martin~R. Bridson and Andr{\'{e}} Haefliger, \emph{Metric spaces of
  non-positive curvature}, Grundlehren der mathematischen Wissenschaften, vol.
  319, Springer, Berlin, 1999.

\bibitem{Cas07}
Christopher~H. Cashen, \emph{Quasi-isometries among tubular groups}, Ph.D.
  thesis, University of Illinois at Chicago, Chicago, IL,
  2007, \href{http://www.math.uic.edu/~ccashen/QIATG.pdf}{\url{http://www.math.u%
ic.edu/~ccashen/QIATG.pdf}}.

\bibitem{CroKle00}
Christopher~B. Croke and Bruce Kleiner, \emph{Spaces with nonpositive curvature
  and their ideal boundaries}, Topology \textbf{39} (2000), 549--556.

\bibitem{MosSagWhy03}
Lee Mosher, Michah Sageev, and Kevin Whyte, \emph{Quasi-actions on trees {I}:
  Bounded valence}, Annals of Mathematics \textbf{158} (2003), no.~1, 115--164,
  \href{http://arXiv.org/abs/math.GR/0010136}{\texttt{arXiv:math.GR/0010136}}.

\bibitem{MosSagWhy04}
\bysame, \emph{Quasi-actions on trees {II}: Finite depth {Bass-Serre} trees},
  preprint, 2004,
  \href{http://arXiv.org/abs/math.GR/0405237}{\texttt{arXiv:math.GR/0405237}}.

\bibitem{PapWhy02}
Panos Papasoglu and Kevin Whyte, \emph{Quasi-isometries between groups with
  infinitely many ends}, Commentarii Mathematici Helvetici \textbf{77} (2002),
  133--144,
  \href{http://arXiv.org/abs/math.GT/0405274}{\texttt{arXiv:math.GT/0405274}}.

\bibitem{ScoWal79}
P.~Scott and C.~T.~C. Wall, \emph{Topological methods in group theory},
  Homological Group Theory (Durham), London Math. Soc. Lecture Notes, September
  1979, pp.~137--203.

\bibitem{Ser03}
Jean-Pierre Serre, \emph{Trees}, Springer Monographs in Mathematics, Springer,
  Berlin, 2003, corrected second printing of the first English edition.

\bibitem{Why99}
Kevin Whyte, \emph{Amenability, bilipschitz equivalence, and the {Von Neumann}
  conjecture}, Duke Mathematical Journal \textbf{99} (1999), no.~1, 93--112,
  \href{http://arXiv.org/abs/math.GR/9704202}{\texttt{arXiv:math.GR/9704202}}.

\bibitem{Why01}
\bysame, \emph{The large scale geometry of the higher {Baumslag-Solitar}
  groups}, Geom. Funct. Anal. \textbf{11} (2001), 1327--1343,
  \href{http://arXiv.org/abs/math.GT/0405272}{\texttt{arXiv:math.GT/0405272}}.

\bibitem{Wis96}
Daniel~T. Wise, \emph{A non-hopfian automatic group}, Journal of Algebra
  \textbf{180} (1996), no.~3, 845--847.

\end{thebibliography}

\end{document}